\def\editmode{0}

\def\bibfilenames{bibman_refs}

\def\spsformat{1}

\if\spsformat1

\documentclass{article}
\usepackage{spconf}
\ninept

\newenvironment{IEEEproof}{\textit{Proof.}}{\hfill $\blacksquare$}

\else
\documentclass[11pt,final,onecolumn]{IEEEtran}
\fi

\usepackage{savesym} %% Needed for compatibility with IEEEtran
\savesymbol{labelindent}
\usepackage{enumitem} % do not include the package enumerate. See below.
\newcommand{\acom}[1]{\textcolor{red}{{[#1]}}} % author comment
\if\editmode1  %%%%%%%%%%%%%%%%%%%%%%%%%%%%%%%%%%%%%%%%%%%%%%%%%%%%%%%%%%%
\usepackage[backend=bibtex,style=alphabetic,sorting=debug,maxbibnames=99]{biblatex}
\DeclareFieldFormat{labelalpha}{\thefield{entrykey}}
\DeclareFieldFormat{extraalpha}{}
\bibliography{\bibfilenames}
\newcommand{\cmt}[1]{\noindent\textcolor{lightgreen}{\underline{[#1]}}} % comment
\newcommand{\hc}[1]{\textcolor{blue}{#1}} % highlight command --> to
\setlistdepth{9}
\newlist{bulletlist}{enumerate}{9}
\setlist[bulletlist,1]{label=$\bullet$}
\setlist[bulletlist,2]{label=$\diamond$}
\setlist[bulletlist,3]{label=$\rightarrow$}
\setlist[bulletlist,4]{label=$\circ$}
\setlist[bulletlist,5]{label=$-$}
\setlist[bulletlist,6]{label=$\square$}
\setlist[bulletlist,7]{label=$\star$}
\setlist[bulletlist,8]{label=$\checkmark$}
\setlist[bulletlist,9]{label=$\Delta$}
\newenvironment{bullets}{\begin{bulletlist}}{\end{bulletlist}}

\newcommand{\blt}[1][noargpassed]{% add a label if an optional argument is passed
  \item%
  \ifthenelse{\equal{#1}{noargpassed}}{}{\cmt{#1}}%
}

\else
\usepackage{cite}
\newcommand{\cmt}[1]{} % comment
\newcommand{\hc}[1]{\textcolor{black}{#1}} % highlight command -->
\newenvironment{bullets}{}{}
\newcommand{\blt}[1][noargpassed]{\ignorespaces}

\fi
\newcommand{\printmybibliography}{
\if\editmode1 
\printbibliography
\else
\bibliography{\bibfilenames}
\fi
}

\usepackage{fixltx2e}

\usepackage{graphicx}

\usepackage{subcaption}              %subfigures with \mbox and \subfigure[]

\usepackage[utf8]{inputenc}

\usepackage{amsfonts}
\usepackage{amsmath}
\usepackage{mathtools}
\usepackage{amssymb}

\usepackage{bm}

\usepackage{color,verbatim}
\usepackage{multirow}
\usepackage{accents}

\usepackage{theoremref}
\usepackage{url}

\newcounter{rulecounter}
\newcommand{\resetrule}{ \setcounter{rulecounter}{0}}
\resetrule

\newtheorem{myauxproblem}{Problem}
\newtheorem{myauxoptionalproblem}{Optional Problem}
\newsavebox{\selvestebox}
\newenvironment{colbox}[1]
  {\newcommand\colboxcolor{#1}%
   \begin{lrbox}{\selvestebox}%
   \begin{minipage}{\dimexpr\columnwidth-2\fboxsep\relax}}
  {\end{minipage}\end{lrbox}%
   \begin{center}
   \colorbox{\colboxcolor}{\usebox{\selvestebox}}
   \end{center}}

\definecolor{orange}{rgb}{1,0.8,0}
\definecolor{gray}{rgb}{.9,0.9,0.9}
\definecolor{darkgray}{rgb}{.3,0.3,0.3}
\definecolor{darkblue}{rgb}{.1,0.0,0.3}
\definecolor{lightblue}{rgb}{0.7,0.7,1}
\definecolor{lightred}{rgb}{1,0.7,.7}
\definecolor{purple}{RGB}{204,153,255}
\definecolor{lightgray}{rgb}{.95,0.95,0.95}
\definecolor{lightgreen}{rgb}{0.3,0.5,0.3}
\definecolor{darkgreen}{rgb}{0.05,0.3,0.05}

\newcommand{\ra}{$\rightarrow$~}

\newcommand{\bbm}[1]{{\bar{\bm #1}}}

\newcommand{\rfield}{\mathbb{R}}

\newcommand{\frob}{{\hc{\text{F}}}}

\newcommand{\sign}{\mathop{\rm sign}}

\newcommand{\transpose}{^\top}
 \newcommand{\define}{:=}
\newcommand{\minimize}{\mathop{\text{minimize}}}

\newcommand{\st}{\mathop{\text{s.t.}}}

\DeclareMathOperator*{\argmin}{arg\,min}

\newcommand{\lagrangian}{{\hc{\mathcal{L}}}}

\newtheorem{myproposition}{Proposition}
\newtheorem{myremark}{Remark}
\newtheorem{myproblemstatement}{Problem Statement}
\newtheorem{mylemma}{Lemma}
\newtheorem{mytheorem}{Theorem}
\newtheorem{mydefinition}{Definition}
\newtheorem{mycorollary}{Corollary}

\usepackage{hyperref} % comment to write a response letter
\usepackage{algorithmic}
\usepackage[ruled,vlined,resetcount,noend,linesnumbered]{algorithm2e}
\usepackage{balance}

\begin{document}

%%%%%%%%%%%%%%%%%%%%%%%%%%%%%%%%%%%%%%%%%%%%%%%%%%%%%%%%%%%%%%%%%%%%% 
\title{Aerial Base Station Placement Leveraging  Radio Tomographic Maps}

\if\spsformat1 \name{ Daniel Romero$^1$, Pham Q. Viet$^1$, and Geert
  Leus$^2$ \thanks{This work was supported by the Research Council of
    Norway through the IKTPLUSS Grant 311994.}}  \address{
  $^1$ University of Agder, Norway. Email \{daniel.romero,viet.q.pham\}@uia.no \\
  $^2$ Delft University of Technology, The Netherlands. Email:
  g.j.t.1eus@tudelft.nl} \else
% \author{Author(s) Name(s)\thanks{Thanks to XYZ agency for funding.}}
\fi

\maketitle
%%%%%%%%%%%%%%%%%%%%%%%%%%%%%%%%%%%%%%%%%%%%%%%%%%%%%%%%%%%%%%%%%%%%% 

\begin{abstract}
  Mobile base stations on board unmanned aerial vehicles (UAVs)
  promise to deliver connectivity to those areas where the
  terrestrial infrastructure is overloaded, damaged, or absent. A
  fundamental problem in this context involves determining a
  minimal set of locations in 3D space where such aerial base
  stations (ABSs) must be deployed to provide coverage to a set of
  users. While nearly all existing approaches rely on average
  characterizations of the propagation medium, this  work
  develops a scheme where the actual channel information is
  exploited by means of a radio tomographic map. A convex
  optimization approach is presented to minimize the number of
  required ABSs while ensuring that the UAVs do not enter no-fly
  regions. A simulation study reveals that the proposed algorithm
  markedly outperforms its competitors.

\end{abstract}

\newcommand{\includefig}[1]{\includegraphics[trim=25 0 50 30,clip,width=0.49\textwidth]{#1}}

\newcommand{\journal}[1]{}
\newcommand{\journalc}[1]{} % journal citation
\newcommand{\conference}[1]{\textcolor{black}{#1}}

\newcommand{\indicator}{\mathbb{I}}
\newcommand{\indicatorinf}{\mathcal{I}}

\newcommand{\groundregion}{{\hc{\mathcal{X}}}}
\newcommand{\height}{{\hc{h}}}
\newcommand{\flyregion}{{\hc{\mathcal{F}}}}

\newcommand{\usernum}{{\hc{M}}}
\newcommand{\userind}{{\hc{m}}}
\newcommand{\absnum}{{\hc{N}}}
\newcommand{\absind}{{\hc{n}}}

\newcommand{\slfgrid}{{\hc{\bar{\mathcal{X}}}}}
\newcommand{\slfgridptind}{{\hc{q}}}
\newcommand{\slfgridptnum}{{\hc{Q}}}
\newcommand{\slfgridside}{{\hc{\slfgridptnum}}_{0}}
\newcommand{\slfgridsidex}{{\hc{\slfgridptnum}}_x}
\newcommand{\slfgridsidey}{{\hc{\slfgridptnum}}_y}
\newcommand{\slfgridsidez}{{\hc{\slfgridptnum}}_z}
\newcommand{\slfgridpt}{{\hc{\bm x}}^\slfgrid}

\newcommand{\grid}{{\hc{\bar{\mathcal{F}}}}}
\newcommand{\gridptind}{{\hc{g}}}
\newcommand{\gridptnum}{{\hc{G}}}    
\newcommand{\gridside}{{\hc{\gridptnum}}_{0}}
\newcommand{\gridsidex}{{\hc{\gridptnum}}_x}
\newcommand{\gridsidey}{{\hc{\gridptnum}}_y}
\newcommand{\gridsidez}{{\hc{\gridptnum}}_z}
\newcommand{\gridpt}{{\hc{\bm x}}^\grid}

\newcommand{\locs}{{\hc{x}}} % loc scalar
\newcommand{\loc}{{\hc{\bm \locs}}}
\newcommand{\userloc}{{\hc{\bm x}}^\text{GT}}
\newcommand{\absloc}{{\hc{\bm x}}^\text{ABS}}
\newcommand{\wavelen}{\hc{\lambda}}
\newcommand{\gain}{\hc{\gamma}}
\newcommand{\shad}{\hc{\xi}}
\newcommand{\slf}{\hc{l}}
\newcommand{\weightfun}{\hc{w}}
\newcommand{\wfwidth}{\hc{\tilde\lambda}} % width of the weight function

\newcommand{\bandwidth}{\hc{W}}
\newcommand{\txpow}{\hc{P_\text{TX}}}
\newcommand{\noisepow}{\hc{\sigma^2}}

\newcommand{\capmat}{\hc{\bm C}}
\newcommand{\capfun}{\hc{C}}
\newcommand{\caps}{\hc{c}}
\newcommand{\caprs}{\hc{\bar c}}
\newcommand{\capvec}{\hc{\bm \caps}}
\newcommand{\caprvec}{\hc{\bbm \caps}} % row vector
\newcommand{\ratemat}{\hc{\bm R}}
\newcommand{\rate}{\hc{ r}}
\newcommand{\raters}{\hc{ \bar r}} % row - scalar
\newcommand{\ratevec}{\hc{\bm \rate}}
\newcommand{\ratervec}{\hc{\bbm \rate}} % row vector
\newcommand{\raterowvec}{\hc{\bbm \rate}}
\newcommand{\ratemin}{\hc{\rate}_\text{min}}
\newcommand{\acti}{\hc{\alpha}}
\newcommand{\actiold}{\hc{\tilde\alpha}}
\newcommand{\activec}{\hc{\bm \alpha}}

\newcommand{\sparsweight}{\hc{ w}}
\newcommand{\sparsweightvec}{\hc{\bm \sparsweight}}

\newcommand{\normvec}{\hc{\bm v}}
\newcommand{\norm}{\hc{ v}}
\newcommand{\vecent}[1]{[#1]}

\newcommand{\slack}{\hc{s}}
\newcommand{\slackvec}{\hc{\bm \slack}}
\newcommand{\slacklow}{\check \slack}
\newcommand{\slackhigh}{\hat \slack}

\newcommand{\Xmat}{{\hc{\bm X}}}
\newcommand{\Zmat}{{\hc{\bm Z}}}
\newcommand{\zrvec}{{\hc{\bbm z}}}
\newcommand{\zs}{{\hc{ z}}} % scalar
\newcommand{\zrs}{{\hc{ \bar z}}} % row - scalar
\newcommand{\zvec}{{\hc{\bm \zs}}}

\newcommand{\Amat}{{\hc{\bm A}}}
\newcommand{\Bmat}{{\hc{\bm B}}}
\newcommand{\Umat}{{\hc{\bm U}}}
\newcommand{\urvec}{{\hc{\bbm u}}}
\newcommand{\us}{{\hc{ u}}}
\newcommand{\urs}{{\hc{ \bar u}}}
\newcommand{\uvec}{{\hc{\bm \us}}}
\newcommand{\Ffun}{{\hc{F}}} % aux function
\newcommand{\Gfun}{{\hc{G}}} % aux function
\newcommand{\admmfunx}{\hc{f}}
\newcommand{\admmfunz}{\hc{h}}

\newcommand{\itnot}[1]{^{#1}}
\newcommand{\itind}{\hc{k}}
\newcommand{\entnot}[1]{[#1]} % entry notation

\newcommand{\gii}{_\gridptind\itnot{\itind}}
\newcommand{\admmstep}{\hc{\rho}}

\newcommand{\lambdas}{{\hc{\lambda}}}
\newcommand{\lambdaslow}{{\hc{\check\lambda}}}
\newcommand{\lambdashigh}{{\hc{\hat\lambda}}}
\newcommand{\nus}{{\hc{\nu}}}
\newcommand{\nuvec}{{\hc{\bm \nus}}}
\newcommand{\mus}{{\hc{ \mu}}}
\newcommand{\muvec}{{\hc{\bm \mus}}}

% \begin{keywords}
%   One, two, three, four, five
% \end{keywords}

\section{Introduction}

\begin{bullets}
  \blt[overview]
  \begin{bullets}
    \blt[motivation]The rapid evolution of the technology of
    unmanned aerial vehicles (UAVs) has spurred extensive research
    to complement terrestrial communication infrastructure
    with base stations mounted on board
    UAVs~\cite{zeng2019accessing}.
    \blt[ABSs] The main use case of such \emph{aerial base
      stations} (ABSs) is to provide connectivity in areas where
    it is insufficient or not available, e.g. because they are
    remote or because of a natural disaster.
    \blt[placement]The research question that arises is at
    which locations one or multiple ABSs need to be deployed to
    provide coverage to the ground terminals (GTs). 
  \end{bullets}

  \blt[literature] 
  \begin{bullets}
    \blt\cmt{one UAV}This question has been extensively
    investigated for a single ABS; see
    e.g. \cite{han2009manet\journalc{,lee2011climbing},boryaliniz2016placement,chen2017map,wang2018adaptive}.    % \begin{bullets}
    %   \blt  \cite{boryaliniz2016placement} mixed integer
    %   program; batch; needs the location of all MUs to start working. 
    %   \blt \cite{chen2017map}. batch. single relay . 
    %   \blt \cite{boryaliniz2019agile} single AirBS.
    %   \blt\cite{han2009manet}
    %   \blt\cite{lee2011climbing}
    % \end{bullets}%
    \blt[multiple UAVs]
    \begin{bullets}%%
      \blt[2D]Other schemes have been proposed to set the 2D
      position of multiple ABSs in a horizontal plane of a given
      height; see
      e.g.~\journal{\cite{lee2010decentralized,andryeyev2016selforganized,galkin2016deployment,lyu2017mounted,romero2019noncooperative,huang2020sparse,mach2021realistic,yin2021multiagent}.}\conference{\cite{romero2019noncooperative}.}
      \blt[3D]In contrast, the focus here is on algorithms capable
      of determining the 3D position of the ABSs. Existing
      works in this context are classified next according to how
      they account for the propagation channel between the ABSs
      and the GTs.
      \begin{bullets}%
        \blt[channel agnostic]\journal{The scheme in
          \cite{park2018formation} could be deemed \emph{channel
            agnostic} since the only way the channel is taken into
          account is in the fact that each ground user associates
          with the ABS from which it receives the strongest
          beacons. This means that it is not possible to know
          whether a spatial arrangement of ABSs is satisfactory
          beforehand, i.e., before the ABSs actually occupy those
          locations, which is clearly inconvenient if one wishes
          to find a close-to-optimal placement. }\conference{First, some
          schemes~\cite{park2018formation} do not model or learn
          the channel and, therefore, the suitability of a
          location cannot be determined before an ABS visits it,
          which drastically increases the time to find a suitable
          placement.}
        \blt[free-space]Besides approaches that assume free-space
        propagation \cite{kim2018topology}, a large number of
        works rely on
        \blt[empirical models\ra based on
        \cite{alhourani2014urban,alhourani2014lap}] the empirical
        model
        from~\cite{alhourani2014urban\journalc{,alhourani2014lap}};
        \begin{bullets}%
          \blt[papers]see e.g.
          \cite{kalantari2016number,hammouti2019mechanism,perabathini2019qos,liu2019deployment,shehzad2021backhaul}.%\cite{mozaffari2016coverage}
          \blt[limitations] The main limitation is that such
          models provide shadowing values in \emph{average
            scenarios}, e.g. in a generic urban environment, but
          are likely to yield highly suboptimal placements in a
          specific environment.
        \end{bullets}%
        \blt[Terrain map/3D model]This limitation is addressed in
        \cite{qiu2020reinforcement,sabzehali2021orientation} by using
        \journal{terrain maps or} 3D models of the deployment
        scenario. \journal{This makes it possible to determine whether
          there is a line of sight between two points without visiting
          them but.}Unfortunately, \journal{terrain maps or} 3D models
        are seldom available and, even when they are, their resolution
        is insufficient for reasonably predicting the channel in conventional
        bands or, for example, when a GT is inside a building. \journal{A typical
          resolution for a terrain map is 30 m, which is therefore
          insufficient in practice. Besides, they cannot accommodate
          the case where a user is inside a building.}
      \end{bullets}%
    \end{bullets}%
  \end{bullets}%

  \blt[contributions]In contrast, the present paper proposes a
  scheme where the air-to-ground channel of the \emph{specific}
  deployment scenario is learned
  \begin{bullets}%
    \blt[tomographic model]
    \begin{bullets}%
      \blt[idea] by relying on the notion of \emph{radio
        tomography}~\cite{patwari2008nesh,patwari2008correlated}.
      \blt[channel map]A \emph{radio map} that provides the
      attenuation between arbitrary points of space is
      constructed based on measurements collected by the GTs and ABSs.
      \journal{and can be updated on-the-go.}%
      \blt[limitations literature]\journal{However, the
        conventional approach to radio tomography entails
        prohibitive complexity for mapping the air-to-ground
        channel, as required here since a 3D grid must be used. In
        the first part of the paper, which deals with constructing
        such radio maps, we show how the cubic complexity of
        this standard approach can be reduced to a linear one
        without sacrificing performance.}\conference{To
        accommodate the special requirements of air-to-ground
        radio maps, the conventional approach to radio
        tomography, which has a cubic complexity in the size of
        the grid, is here replaced with a linear complexity
        algorithm.}
    \end{bullets}%
    \blt[sparse optimization]\journal{The second part of the paper
      assumes that the radio map is given and finds the optimal
      ABS locations according to a convex relaxed version of a
      criterion that aims at minimizing the number of ABSs. This
      approach adopts a 3D discretization of the airspace, which
      makes it possible to accommodate flight restrictions such as
      those imposed by obstacles (e.g. buildings) and no-fly
      zones. This is a key advantage over most of the algorithms
      in the literature, which cannot accommodate such
      constraints.}\conference{Using this radio map, a placement
      algorithm is proposed to minimize the number of ABSs
      required to guarantee a minimum rate for all GTs. Unlike
      most competing algorithms, it is based on a convex program,
      it can accommodate no-fly zones, and has low computational
      complexity. }%
    \blt[simulator]The third contribution is an open source
    simulator\footnote{\scriptsize{\url{https://github.com/uiano/abs_placement_via_radio_maps}} }
    that allows testing and developing  algorithms for ABS placement; see
    Fig.~\ref{fig:urban}.

  \end{bullets}%

  \begin{figure}[!t]
    \centering
    \includegraphics[width=0.48\textwidth]{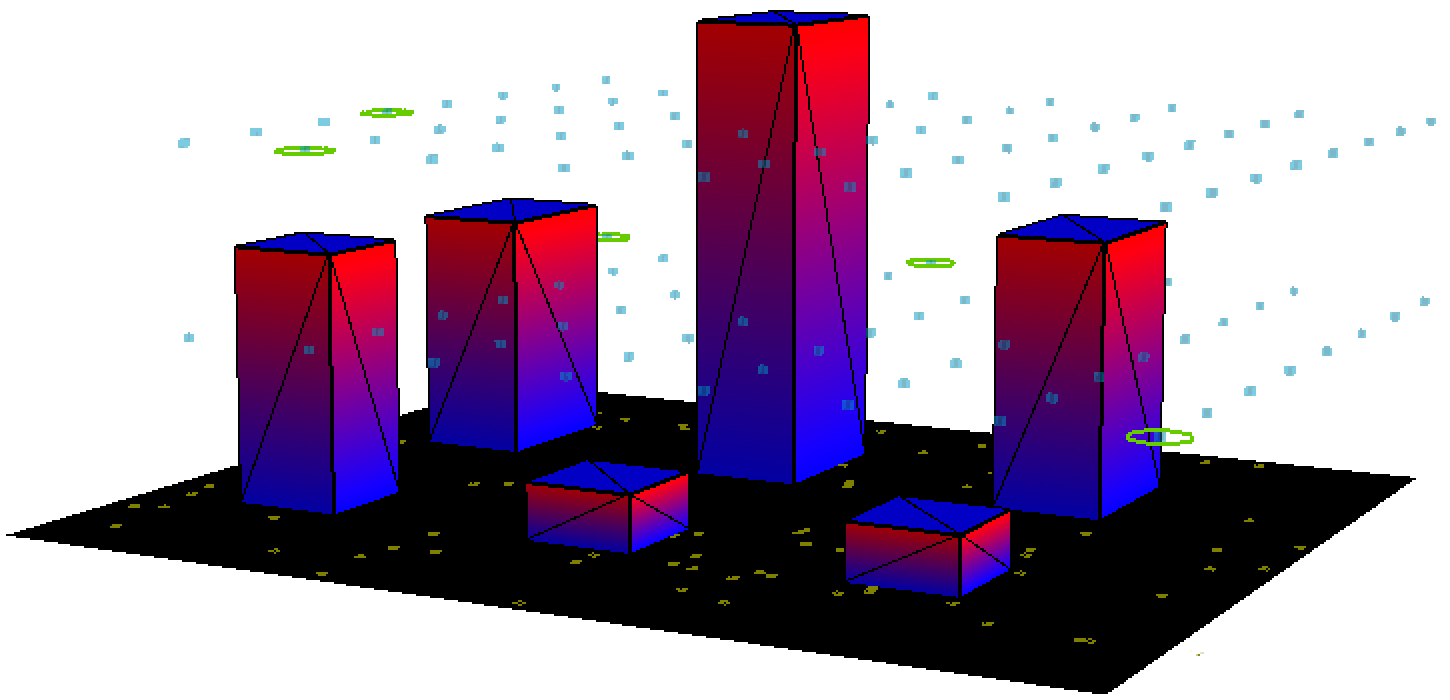}
    \caption{Example of ABS placement in an urban environment with the
      developed simulator. GTs are represented by markers on the
      ground, grid points by blue dots, and ABS positions by green
      circles.}
    \label{fig:urban}
  \end{figure}

  \blt[paper structure]\emph{Paper structure.}
  Sec.~\ref{sec:model} \journal{introduces the model }and formulates the
  problem. The construction and evaluation of radio maps is described
  in Sec.~\ref{sec:radiomaps}. An algorithm for ABS placement using
  radio maps is then proposed in Sec.~\ref{sec:placement}. Performance
  evaluation is carried out in Sec.~\ref{sec:exp} by means of the
  developed simulator. Finally, Secs.~\ref{sec:related}
  and~\ref{sec:conclusions} respectively discuss the related work and
  present the main conclusions. The supplementary material contains an
  algorithm for approximating tomographic integrals and the derivation
  of the placement algorithm.
  
  \blt[notation]\emph{Notation.}
  \begin{bullets}%
    \blt[] $\rfield_+$ is set of non-negative real numbers.
    \blt[] Boldface uppercase (lowercase) letters denote matrices
    (column vectors).
    \blt[] $a\entnot{i}$ represents the $i$-th entry of vector $\bm a$.
    \blt[] Notation $\bm 0$ (respectively $\bm 1$) refers to the
    matrix of  the appropriate dimensions with all zeros (ones).
    \blt[] $\|\bm A\|_\frob$ denotes Frobenius norm of matrix $\bm A$,
    whereas $\|\bm a\|_p$ denotes the $\ell_p$-norm of vector $\bm
    a$. With no subscript, $\|\bm a\|$ stands for the
    $\ell_2$-norm. Inequalities between vectors or matrices must be
    understood entrywise.

  \end{bullets}%
\end{bullets}%

  \section{Model and Problem Formulation}
  \label{sec:model}

\begin{bullets}%
  \blt[model]
  \begin{bullets}%
    \blt[users]Consider $\usernum$ users or ground terminals (GTs)
    located at positions
    $\{\userloc_1, \ldots,\userloc_\usernum\} \subset
    \groundregion \subset\rfield^3$, where region $\groundregion$
    will typically include points on the ground and inside
    buildings.
    \blt[abss]To provide connectivity to the GTs, $\absnum$ ABSs
    are deployed at positions
    $\{\absloc_1,\ldots,\absloc_\absnum \}\subset\flyregion\subset
    \rfield^3$, where  $\flyregion$ comprises all
    locations where a UAV is allowed to fly. This excludes no-fly
    zones, airspace occupied by buildings, and altitudes out
    of legal limits.
    \blt[channel]%
    \begin{bullets}%%
      \blt[gain]To simplify the exposition, the focus will be on the
      downlink and it will be assumed that the channel is not frequency
      dispersive. The rate of the communication link between the
      $\userind$-th GT and an ABS at position
      $\absloc\in\groundregion$ is determined by the channel gain
      and noise power. The former is given by
      \begin{align}
        \label{eq:gain}
        \gain_\userind(\absloc) = 20\log_{10}\left(
        \frac{\wavelen}{4\pi \|\userloc_\userind- \absloc\|}
        \right) - \shad(\userloc_\userind, \absloc),
      \end{align}
      where $\wavelen$ is the wavelength associated with the
        carrier frequency of the transmission and function $\shad$ denotes shadowing. Small-scale
      fading is  ignored for simplicity, but the ensuing
        formulation can be adapted to accommodate the associated
        uncertainty.  \blt[capacity]The capacity~is
      \journal{\begin{align}
                 \label{eq:capfun}
                 \capfun_\userind(\absloc) = \bandwidth \log_{2}\left(
                 1 + \frac{\txpow 10^{\gain_\userind( \absloc)/10}}{\noisepow}
                 \right) , 
               \end{align}}
             \conference{
               \begin{align}
                 \label{eq:capfun}
                 \capfun_\userind(\absloc) = \bandwidth \log_{2}\left(
                 1 + {\txpow 10^{\gain_\userind( \absloc)/10}}/{\noisepow}
                 \right) , 
               \end{align}
             }where $ \bandwidth$ denotes bandwidth,
             \journal{$\txpow$ the (frequency flat) transmit power
               spectral density (PSD), and $\noisepow$ the noise
               PSD.}\conference{$\txpow$ the transmit power, and
               $\noisepow$ the noise power.}
             \blt[comb. rate]Since the $\userind$-th GT may connect to
             one or multiple ABSs,  it may receive
             a rate up to
             $\sum_\absind\capfun_\userind(\absloc_\absind)$.
             \blt[unlimited backhaul]As usual in the literature, it is
             assumed that the backhaul connection of the ABSs
             has sufficiently high capacity, yet the proposed scheme can
             be generalized to accommodate backhaul constraints.

             % 
             % \blt[rate]As usual in cellular networks, $\txpow$ is
             % constant but $\bandwidth$ is adjusted to allocate the
             % desired rate to each user. To guarantee a minimum rate
             % $\ratemin$, multiple ABSs may provide connectivity to a
             % given user. With $\{\absloc_\absind\}_{\absind=1}^\absnum$
             % indicating the locations of the ABSs, the rates
             % $\rate(\userloc, \absloc_\absind)$ are chosen so that
             % \begin{align}
             %   \sum_\absind \rate(\userloc, \absloc_\absind)
             %   \geq \ratemin
             % \end{align}
             % while
             % $0\leq \rate(\userloc, \absloc_\absind)\leq
             % \capfun(\userloc, \absloc_\absind)$ provided that the
             % $\{\absloc_\absind\}_{\absind=1}^\absnum$ are such that this is possible.          

           \end{bullets}%
         \end{bullets}%

         \blt[problem formulation]
         \begin{bullets}%
           \blt[criterion]The problem is to find a minimal set of ABS
           locations that guarantees a minimum rate for every
           user. This criterion arises naturally in some of the main
           use cases of UAV-assisted networks such as emergency
           response or disaster management.  \blt[problem]Formally,
           the problem can be stated as follows:
           \begin{subequations}
             \label{eq:problemf}
             \begin{align}
               \minimize_{\absnum,
               \{\absloc_\absind\}_{\absind=1}^\absnum}~ &\absnum\\
               \label{eq:problemfcap}
               \st~& \textstyle\sum_\absind \capfun_\userind(\absloc_\absind)
                         \geq \ratemin,~ \userind=1,\ldots,\usernum,    \\
               \label{eq:problemfpos}
                                                             & \absloc_\absind \in \flyregion,~~  \absind=1,\ldots,\absnum.
             \end{align}
           \end{subequations}
           To simplify notation, the same rate $\ratemin$ is assumed
             across GTs, but different rates can be set up to
             straightforward modifications.
         \end{bullets}%
       \end{bullets}%

       \section{Tomographic Radio Maps }
       \label{sec:radiomaps}

\begin{figure}[!t]
         \centering
         \includegraphics[width=.49\textwidth]{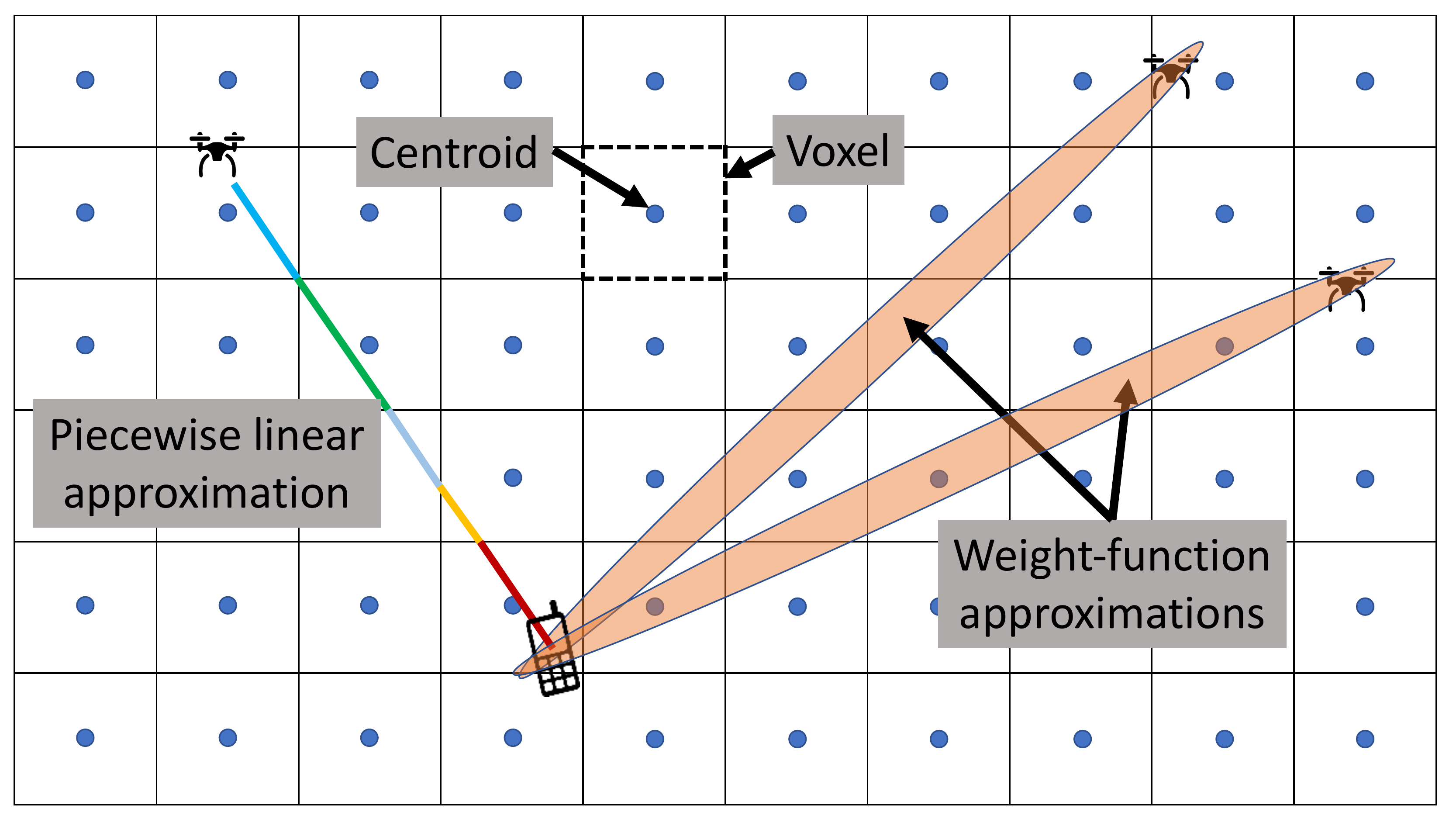}
         \caption{2D illustration of the conventional weight-function
           approximation of the tomographic integral
           \eqref{eq:tomoint} (orange ellipses) and the approximation
           adopted here (colored line segment). Observe that the upper
           ellipse contains no centroid and, therefore, the
           approximation will yield zero attenuation regardless of the
           values of the SLF. }
         \label{fig:tomography}
       \end{figure}

       \begin{bullets}%

         \blt[overview] 
         \begin{bullets}%
           \blt[difficulty\ra unk. shad.]The first difficulty when solving
           \eqref{eq:problemf} is that the function
           $\capfun_\userind(\absloc)$ is unknown since the shadowing term
           $\shad(\userloc_\userind, \absloc)$ in \eqref{eq:gain} is unknown.
           \blt[proposed]The approach proposed here is to rely on a radio map
           that provides $\shad(\userloc, \absloc)$ for all $\userloc$ and
           $\absloc$.
           \blt[tomogr.\ra ch. maps]Such a map can be constructed by means of
           the so-called tomographic (or NeSh) model~\cite{patwari2008nesh},
           as considered in the literature of channel-gain cartography; see
           \cite{romero2018blind} and references therein. However, the
           existing works in this context focus on ground-to-ground
           channels. Constructing radio maps of air-to-ground channels
           involves special challenges that render existing approaches
           unsuitable, as discussed later.
         \end{bullets}%

         \blt[NeSh model]The radio tomographic model~\cite{patwari2008nesh}
         prescribes that 
         \begin{bullets}%
           \blt[integral]
           \begin{align}
             \label{eq:tomoint}
             \shad(\loc_1,\loc_2) = \frac{1}{
             \| \loc_1 -\loc_2\|_2^{1/2}
             } \int_{\loc_1}^{\loc_2}\slf(\loc)d\loc,
           \end{align}
           where  the function $\slf$ inside the line integral is termed \emph{spatial
             loss field} (SLF) and quantifies the local attenuation (absorption) that
           a signal suffers at each position. 
           \blt[estimation]The SLF can be estimated in a first stage before
           solving~\eqref{eq:problemf} by collecting measurements of the form
           $(\absloc, \userloc, \gain_\userind(\absloc))$ and applying
           standard estimation techniques; see
           e.g.~\cite{wilson2009regularization,kanso2009compressed,romero2018blind}.

           \blt[evaluation]In practice, to estimate $\slf$ and evaluate
           \eqref{eq:tomoint}, function $\slf$ needs to be discretized by
           storing its values
           $\slf(\slfgridpt_1),\ldots,\slf(\slfgridpt_\slfgridptnum)$ on a 3D regular
           grid of $\slfgridptnum$ points
           $\slfgrid\define\{\slfgridpt_1,\ldots,\slfgridpt_\slfgridptnum\}$.
           \begin{bullets}%
             \blt[conventional]
             \begin{bullets}%
               \blt[description]The conventional approach approximates
               \eqref{eq:tomoint} as a weighted
               sum~\cite{hamilton2014modeling} of the values
               $\slf(\slfgridpt_\slfgridptind)$ for which the centroid
               $\slfgridpt_\slfgridptind$ lies inside an ellipsoid
               with foci at $\loc_1$ and $\loc_2$; see the ellipses in
               Fig.~\ref{fig:tomography} for a depiction in 2D.
               \blt[limitations]
               \begin{bullets}%
                 \blt[discontinuous] Unfortunately, it can be easily seen
                 from Fig.~\ref{fig:tomography} that the resulting approximation
                 of $\shad(\loc_1,\loc_2)$ is a discontinuous function of
                 $\loc_1$ and $\loc_2$. It may even be 0 even when
                 $\slf(\slfgridpt_\slfgridptind)\neq0~\forall \slfgridptind$.
                 \blt[num gridpts]To minimize these effects, the grid point
                 spacing needs to be small relative to the length of the
                 minor axis, which is commonly set in the order of the
                 wavelength. Thus, for standard centimetric wavelengths and
                 regions $\groundregion$ with sides in the order of km and height
                 in the order of 100 m, $\slfgridptnum$ must be in the order
                 of $10^{14}$, which is prohibitively high.
                 \blt[complexity]Finally, the complexity of such an
                 approximation is $\mathcal{O}(\slfgridside^3)$ for a
                 $\slfgridside\times\slfgridside\times\slfgridside$ grid. 
               \end{bullets}%
             \end{bullets}%
             
             \blt[proposed]To remedy these issues, this paper advocates
             approximating the integral in \eqref{eq:tomoint} as a line
             integral of a piecewise constant approximation of $\slf$, as
             already hinted in~\cite{kanso2009compressed} for tomographic
             imaging. This involves obtaining the intersections between the
             the voxel boundaries and the line segment that connects the
             transmitter to the receiver locations; see the colored
             segment in
             Fig.~\ref{fig:tomography}. A possible implementation along the
             lines of \cite[Sec. I-B-1]{mitchell1990comparison} is presented
             in the supplementary material, but others are possible. The
             resulting approximation is continuous, can be used with large
             grid point spacing, and can be computed with complexity only
             $\mathcal{O}(\slfgridside)$ for a
             $\slfgridside\times\slfgridside\times\slfgridside$
             grid. \journal{\acom{draw connections with Venkat and
                 Leus' work on rainfall estimation using microwave links}}

           \end{bullets}%
         \end{bullets}%
       \end{bullets}%

       \section{Placement with Min-rate Guarantees}
       \label{sec:placement}

       % \begin{subequations}
       %   \begin{align}
       %     \minimize_{\activec\in \{ 0, 1\}^{ \gridptnum}}\quad & \sum_{\gridptind=1}^{\gridptnum}\acti_\gridptind\\
       %     \st \quad & \sum_{\gridptind=1}^{\gridptnum}\acti_\gridptind\capvec_\gridptind \geq \ratemin \bm 1 
       %   \end{align}
       % \end{subequations}

       % \subsection{Solution via ADMM}

       \begin{bullets}%

         \blt[overview]

         \begin{bullets}%
           \blt[channel known]The approach in Sec.~\ref{sec:radiomaps} makes it
           possible to find the shadowing between any two points and,
           therefore, the channel gain and capacity; cf.~\eqref{eq:gain} and
           \eqref{eq:capfun}. The constraint in
           \eqref{eq:problemfcap} can thus be evaluated.
           \blt[nonconvex]Yet, solving \eqref{eq:problemf} is
             challenging: even if $\absnum$ were known and one just
             needed to find feasible
             $\{\absloc_\absind\}_{\absind=1}^\absnum$, the problem would
             still be non-convex due to the constraints.
           \blt[discretization]To bypass this difficulty, the proposed approach
           involves discretizing the flight region $\flyregion$ into a \emph{flight
           grid}
           $\grid\define \{\gridpt_1,\ldots,\gridpt_\gridptnum \}\subset\flyregion\subset
           \rfield^3$; see~Fig.~\ref{fig:urban}. % For simplicity, the same
           % notation is used  for the flight grid and the SLF grid in
           % Sec.~\ref{sec:radiomaps}, but these grids do not generally coincide.  
         \end{bullets}%

         \blt[reformulation]
         \begin{bullets}%
           \blt[constraint]Since $\grid$ contains only points where an ABS
           can be placed, solving \eqref{eq:problemf} amounts to finding the
           smallest subset of $\absnum$ points of $\grid$ that satisfies
           \eqref{eq:problemfcap}. To see this, replace
           $\absloc_\absind \in \flyregion$ in \eqref{eq:problemfpos} with
           $\absloc_\absind \in \grid$ and let $\acti_\gridptind$ be 1
           if there is an
           ABE at $\gridpt_\gridptind$ and 0 otherwise. The summation in
           \eqref{eq:problemfcap} can then be expressed as
           $\journal{\sum_\absind \capfun_\userind(\absloc_\absind) =}\sum_\gridptind
           \acti_\gridptind\capfun_\userind(\gridpt_\gridptind)$. Since the
           number of ABSs can be written as
           $\sum_\gridptind \acti_\gridptind$, the discretized version of
           \eqref{eq:problemf} becomes \journal{
             \begin{subequations}
               \label{eq:probactisum}
               \begin{align}
                 \minimize_{\activec\in \{ 0, 1\}^{ \gridptnum}}\quad & \sum_{\gridptind=1}^{\gridptnum}\acti_\gridptind\\
                 \st \quad &
                             \sum_{\gridptind=1}^{\gridptnum}\acti_\gridptind\capvec_\gridptind \geq \ratemin \bm 1,
               \end{align}
             \end{subequations}}
           \conference{
             \begin{subequations}
               \label{eq:probactisum}
               \begin{align}
                 \minimize_{\activec\in \{ 0, 1\}^{ \gridptnum}}\quad & \textstyle\sum_{\gridptind}\acti_\gridptind\\
                 \st \quad &
                             \textstyle\sum_{\gridptind}\acti_\gridptind\capvec_\gridptind \geq \ratemin \bm 1,
               \end{align}
             \end{subequations}    
           }where
           $\capvec_\gridptind\define[\capfun_1(\gridpt_\gridptind),\ldots,
           \capfun_\usernum(\gridpt_\gridptind)]\transpose$.
         \end{bullets}%%
         \blt[approaches]%
         \begin{bullets}%%
           \blt[combinatorial]Problem \eqref{eq:probactisum} is of a combinatorial
           nature and can be solved for small $\gridptnum$ by exhaustive
           search. However, the complexity of such a task is exponential
           and, therefore, it is preferable to adopt an approximation that
           can be efficiently computed.
           \blt[interior]One possibility is to relax the constraint
           $\activec\in \{ 0, 1\}^{ \gridptnum}$ as well as the objective
           and apply an interior-point solver. This approach is described
           in the supplementary material
           \blt[why not interior]but not pursued here due to the well-known poor
           scalability of this kind of methods with the number of
             variables and constraints~\cite{lin2021admm}. Indeed, in this application,
           $\gridptnum$ can be in the order of millions, which would
           render the cubic complexity of interior-point methods
           prohibitive. 
           \blt[]Instead, this section presents a solver based on the
           \emph{alternating-direction method of multipliers}
           (ADMM)~\cite{boyd2011distributed} whose complexity is linear
           in~$\gridptnum$.
         \end{bullets}%

         \blt[reformulation]
         \begin{bullets}%
           \blt[chg variables]
           Suppose that there exists no grid point such
           that $\capvec_\gridptind=\bm 0$. Otherwise, $\gridpt_\gridptind$ can be
           disregarded without further implications. By applying the
           change of variables
           $\acti_\gridptind\capvec_\gridptind \rightarrow
           \ratevec_\gridptind$, it is clear that Problem
           \eqref{eq:probactisum} can be equivalently written as
           \journal{
             \begin{subequations}
               \label{eq:sumindis}
               \begin{align}
                 \minimize_{\ratemat\in \rfield^{\usernum \times \gridptnum}}\quad & \sum_{\gridptind=1}^{\gridptnum}  \indicator[\ratevec_\gridptind\neq \bm 0]\\
                 \st \quad & \sum_{\gridptind=1}^{\gridptnum}\ratevec_\gridptind \geq \ratemin \bm 1 \\
                                                                                   & \ratevec_\gridptind \in \{\bm 0, \capvec_\gridptind\},~\gridptind=1,\ldots,\gridptnum,
               \end{align}
             \end{subequations}}
           \conference{
             \begin{subequations}
               \label{eq:sumindis}
               \begin{align}
                 \minimize_{\ratemat\in \rfield^{\usernum \times \gridptnum}}\quad & \textstyle\sum_{\gridptind=1}^{\gridptnum}  \indicator[\ratevec_\gridptind\neq \bm 0]\\
                 \label{eq:sumindisratemin}
                 \st \quad &
                             \textstyle\sum_{\gridptind=1}^{\gridptnum}\ratevec_\gridptind \geq \ratemin \bm 1 \\
                 \label{eq:sumindiscapacity}
                                                                                   & \ratevec_\gridptind \in \{\bm 0, \capvec_\gridptind\},~\gridptind=1,\ldots,\gridptnum,
               \end{align}
             \end{subequations}
           }where     $\ratemat\define[\ratevec_1,\ldots,\ratevec_\gridptnum]$ and
           $\indicator[\cdot]$ is a function that returns 1 when the
           condition in brackets holds and 0 otherwise.
           \blt[relaxing the constraint] It will now be argued that
           relaxing the constraint
           $\ratevec_\gridptind \in \{\bm 0, \capvec_\gridptind\}$ as
           $\bm 0 \leq \ratevec_\gridptind \leq \capvec_\gridptind$
           entails no loss of optimality. On the one hand, if
           $\{\ratevec_\gridptind\}_\gridptind$ are feasible for
           \eqref{eq:sumindis}, then they are feasible for the relaxed
           problem and yield the same objective value.  On the other
           hand, if $\{\ratevec_\gridptind\}_\gridptind$ are feasible
           for the relaxed problem, setting those non-zero
           $\ratevec_\gridptind$ equal to $\capvec_\gridptind$ yields
           a feasible point for \eqref{eq:sumindis} that attains the
           same objective value.

           \blt[min rate constr] The next step is to show that, after
           relaxing \eqref{eq:sumindiscapacity}, the inequality in
           \eqref{eq:sumindisratemin} can be replaced with an equality
           without loss of optimality. First, note that
           \eqref{eq:sumindisratemin} can be written as
           $\ratemat \bm 1 \geq \ratemin \bm 1$. Upon letting
           $\raterowvec_\userind\in \rfield^\gridptnum$ denote the
           $\userind$-th column of $\ratemat\transpose$, % , i.e.,
           % $\ratemat\transpose = [\raterowvec_1,\ldots,
           % \raterowvec_\usernum]$
           constraint \eqref{eq:sumindisratemin} becomes
           $\raterowvec_\userind\transpose\bm 1 \geq
           \ratemin,~\userind=1,\ldots,\usernum$.  Now consider a
           feasible $\ratemat$ and note that if
           $\raterowvec_{\userind_0}\transpose\bm 1 > \ratemin$ for
           some $\userind_0$, then replacing
           $\raterowvec_{\userind_0}$ with
           $\raterowvec_{\userind_0}'\define \ratemin
           \raterowvec_{\userind_0}/ (\bm 1\transpose
           \raterowvec_{\userind_0})$ yields another feasible
           $\ratemat'$ that satisfies
           $(\raterowvec_{\userind_0}')\transpose\bm 1 = \ratemin$ and
           that attains the same objective value as
           $\ratemat$. Applying this logic for all $\userind$ yields a
           feasible matrix that satisfies $\ratemat \bm 1 = \ratemin \bm 1$
           without affecting the objective value.

\newcommand{\widemargin}[1]{\parbox{\dimexpr\textwidth-2\algomargin\relax}{#1}}

\begin{algorithm}[t]
  \SetAlgoLined
  \KwData{$\capmat\in \rfield_+^{\usernum \times \gridptnum}$,
    $\ratemin\in \rfield_+$, $\{\sparsweight_\gridptind\}_\gridptind\subset
    \rfield_+$,~$\admmstep>0$}
  \textbf{Initialize}~ $\Umat\itnot{1}\in \rfield_+^{\usernum \times \gridptnum}$ and $\Zmat\itnot{1}\in \rfield_+^{\usernum \times \gridptnum}$ \\
  \For{$\itind=1,2,\ldots$}{
    \For{$\gridptind=1,2,\ldots,\gridptnum$}{
      Bisection: find $\slack_\gridptind\itnot{\itind+1}$ s.t. \widemargin{$\bm 1\transpose \max(\zvec_\gridptind\itnot{\itind
      }-\uvec_\gridptind\itnot{\itind} -\slack_\gridptind\itnot{\itind+1}\bm 1, \bm 0) = {\sparsweight_\gridptind}/{\admmstep}$}\\
      Set $\ratevec_\gridptind\itnot{\itind+1} = \min(\zvec_\gridptind\itnot{\itind
      }-\uvec_\gridptind\itnot{\itind}, \slack_\gridptind\itnot{\itind+1}\bm 1)$\\
    }
    \For{$\userind=1,2,\ldots,\usernum$}{
      Bisection: find $\lambdas$~s.t.
      \widemargin{$\bm 1\transpose\max(\bm 0,    \min(\caprvec_\userind, \ratervec_\userind\itnot{\itind+1} +
        \urvec_\userind\itnot{\itind} - \lambdas\bm 1)) =\ratemin $}\\[1ex]
      \widemargin{Set $    \zrvec_\userind\itnot{\itind+1} =
\max(\bm 0,    \min(\caprvec_\userind, \ratervec_\userind\itnot{\itind+1} +
    \urvec_\userind\itnot{\itind} - \lambdas\bm 1))
$}\\
}
Set $\Umat\itnot{\itind+1} = \Umat\itnot{\itind} +
    \ratemat\itnot{\itind+1}  -
    \Zmat\itnot{\itind+1}$\\
    \textbf{If} convergence(~) \textbf{then}  return $\ratemat\itnot{\itind+1}$
    % read current\;
    % \eIf{understand}{
    %   go to next section\;
    %   current section becomes this one\;
    % }{
    %   go back to the beginning of current section\;
    % }
    % \lIf{condition}{\Return{} \True}
  }
  \caption{ABS Placement}
\label{algo:placement}    
\end{algorithm}

           \blt[objective]The objective
           $\sum_{\gridptind=1}^{\gridptnum}
           \indicator[\ratevec_\gridptind\neq \bm 0]$ can be equivalently
           expressed as
           $\sum_{\gridptind=1}^{\gridptnum}
           \indicator[\|\ratevec_\gridptind\|_\infty\neq \bm 0]$, where the
             $\ell_\infty$-norm $\|\bm v\|_\infty$ equals the largest absolute
             value of the entries of vector $\bm v$. Clearly,
           $\sum_{\gridptind=1}^{\gridptnum}
           \indicator[\|\ratevec_\gridptind\|_\infty\neq \bm 0] = \|
           [\|\ratevec_1\|_\infty,\ldots,
           \|\ratevec_\gridptnum\|_\infty]\transpose\|_0$, which suggests the
           relaxation
           $ \| [\|\ratevec_1\|_\infty,\ldots,
           \|\ratevec_\gridptnum\|_\infty]\transpose\|_1=\sum_\gridptind
           \|\ratevec_\gridptind\|_\infty$, or its reweighted version
           $\sum_\gridptind \sparsweight_\gridptind
           \|\ratevec_\gridptind\|_\infty$, where
           $\{\sparsweight_\gridptind\}_\gridptind$ are non-negative
           constants set as in~\cite{candes2008reweighted}.
           \blt[resulting prob.]With these observations, the problem becomes
           \begin{subequations}
             \label{eq:ubeqproblem}
             \begin{align}
               \minimize_{\ratemat\in \rfield^{\usernum \times \gridptnum}}\quad & 
                                                                                   \textstyle    \sum_\gridptind \sparsweight_\gridptind
                                                                                   \|\ratevec_\gridptind\|_\infty \\
               \st \quad & \ratemat \bm 1 = \ratemin \bm 1,~ \bm 0 \leq \ratemat \leq \capmat,
             \end{align}
           \end{subequations}
           \blt[interpretation]where the $(\userind,\gridptind)$-th
           entry of  $\capmat\in \rfield_+^{\usernum\times\gridptnum}$ is
           given by
           $\caps_{\userind,\gridptind}\define\capfun_\userind(\gridpt_\gridptind)$,
           i.e., the capacity of the link between the $\userind$-th user and
           the $\gridptind$-th grid point. The $(\userind,\gridptind)$-th entry
           of $\ratemat$ therefore satisfies
           $0\leq \rate_{\userind,\gridptind}\leq \caps_{\userind,\gridptind}$,
           which means that it can be interpreted as the rate at which a
           \emph{virtual ABS} placed at grid point $\gridpt_\gridptind$
           communicates with the $\userind$-th user. In case that
           $\rate_{\userind,\gridptind}=0$ for all $\userind$, then no \emph{actual}
           ABS needs to be deployed at $\gridpt_\gridptind$. In other words,
           the virtual ABS at  $\gridpt_\gridptind$ corresponds to an actual
           ABS only if $\rate_{\userind,\gridptind}\neq 0$ for some
           $\userind$. 

           % In other words,
           % the number $\absnum$ of ABSs equals the number of grid points
           % $\gridpt_\gridptind$ for which $\rate_{\userind,\gridptind}\neq 0$
           % for some $\userind$ or, equivalently, the number of grid points
           % $\gridpt_\gridptind$ for which
           % $\ratevec_\gridptind\define
           % [\rate_{1,\gridptind},\ldots,\rate_{\usernum,\gridptind}]\transpose\neq
           % \bm 0$. This number can be expressed as
           % $\sum_{\gridptind=1}^{\gridptnum} \indicator[\ratevec_\gridptind\neq
           % \bm 0]$.

           % \blt[constraint]The constraint \eqref{eq:problemfcap} can be
           % expressed with this notation as
           % $\sum_\gridptind \rate_{\userind,\gridptind}\geq \ratemin~\forall
           % \userind$ or, in matrix form as
           % $\ratemat \bm 1 \geq \ratemin \bm 1$, where
           % $\ratemat\define[\ratevec_1,\ldots,\ratevec_\gridptnum]$. The
           % resulting discretized version of \eqref{eq:problemf} reads as
           % \begin{subequations}
           %   \label{eq:ubproblem}
           %   \begin{align}
           %     \minimize_{\ratemat\in \rfield^{\usernum \times \gridptnum}}\quad & \sum_{\gridptind=1}^{\gridptnum}  \indicator[\ratevec_\gridptind\neq \bm 0]\\
           %     \st \quad & \ratemat \bm 1 \geq \ratemin \bm 1, \\
           %     & \bm 0 \leq \ratemat \leq \capmat.
           %   \end{align}
           % \end{subequations}

         \end{bullets}%
         \blt[]
       \end{bullets}%  

       Within the ADMM framework, Problem~\eqref{eq:ubeqproblem} can
       be decomposed into one subproblem per row and column of
       $\ratemat$. Each problem involves solving a bisection task of a
       1D monotonically decreasing function and therefore can be solved with
       $\mathcal{O}(1)$ evaluations. The total complexity is
       $\mathcal{O}(\usernum\gridptnum)$, much smaller than the
       $\mathcal{O}((\gridptnum+2\usernum)^3)$ complexity per inner
       iteration of an interior-point method; cf. the supplementary
       material. The algorithm is shown as Algorithm~\ref{algo:placement} and
       is derived in the supplementary material. In the notation used
       therein, if $\bm A$ is a matrix, then $\bm a_m$ is its $m$-th
       column and $\bbm a_n\transpose$ its $n$-th row. Furthermore,
       superscripts indicate the iteration index, $\admmstep>0$ is
       the step size, and the $\min$ and $\max$ operators act
       entrywise.

       % It is well known that interior-point methods are accurate but their
       % scalability is limited. For that reason, in contrast ADMM algorithms
       % offer a greater scalability, which is necessary when the number of
       % grid points is large, as will typically be case in ABS placement
       % applications.

       \section{Numerical Experiments}
       \label{sec:exp}

       \begin{figure}[!t]
         \centering
         \includefig{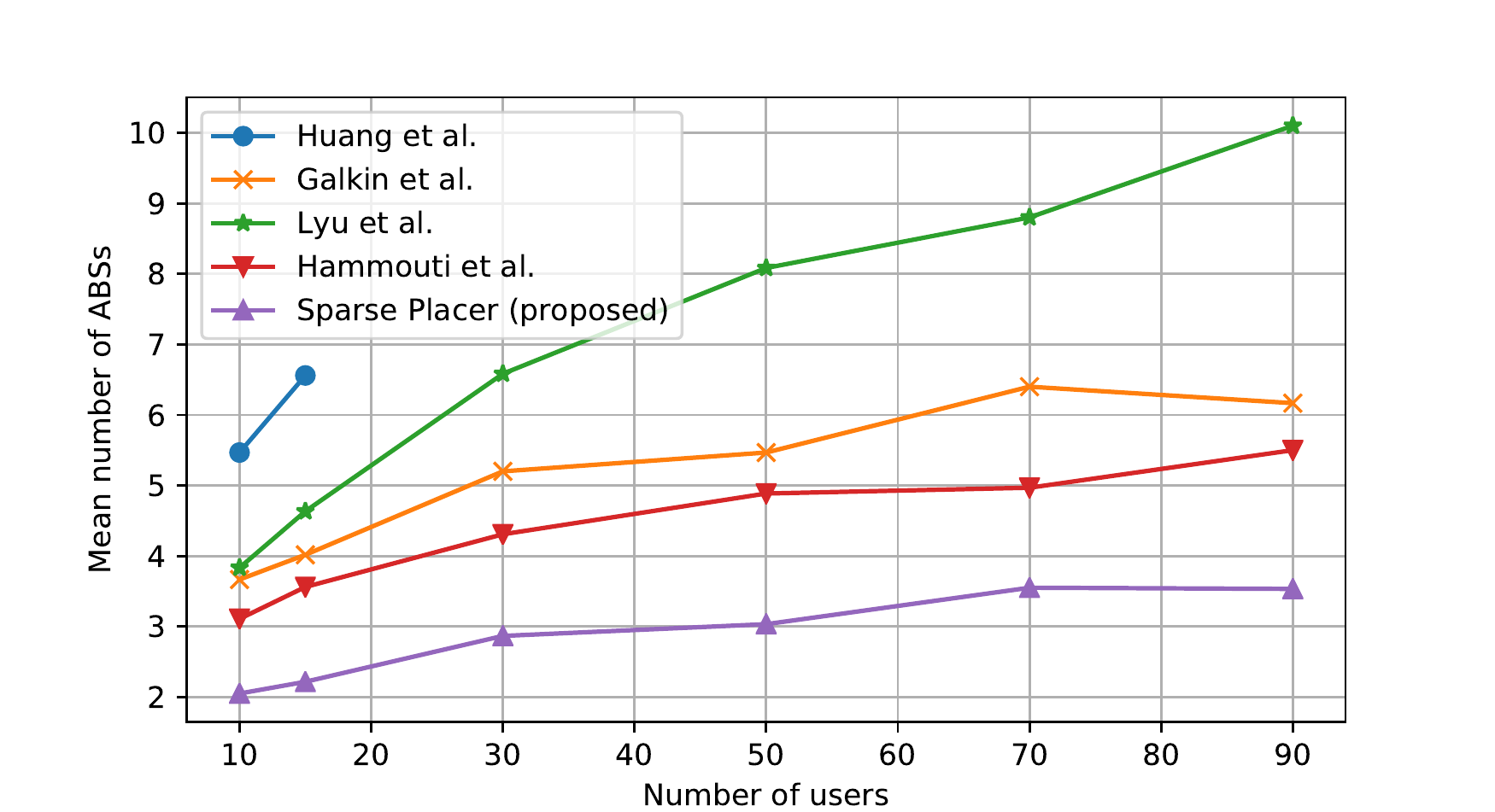}
         \caption{Mean minimum number of ABSs required to provide a minimum
           rate of  $\ratemin=5$ Mb/s vs. the number of GTs ($\height=53$ m, 
           $20\times 30\times 5$ SLF grid, $9\times 9\times 3$ fly grid).
         }
         \label{fig:num_users}
       \end{figure}

       \begin{figure}[!t]
         \centering
         \includefig{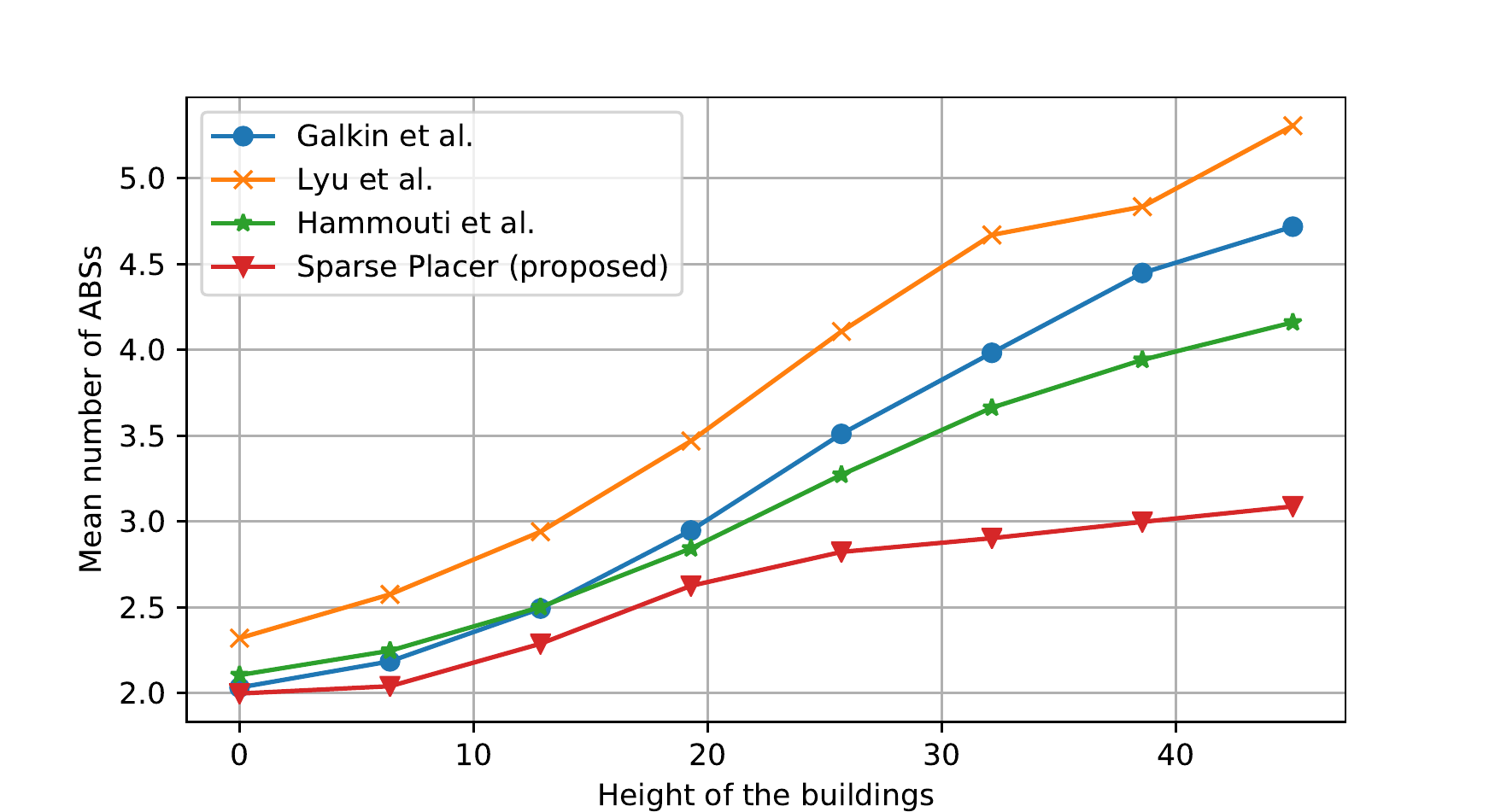}
         \caption{Mean minimum number of ABSs required to provide a minimum
           rate of   $\ratemin=20$ Mb/s vs. $\height$ [m] ($20\times 30\times 150$ SLF
           grid, $9\times 9\times 5$ fly grid).}
         \label{fig:height}
       \end{figure}

       \cmt{simulation setup}
       \begin{bullets}%
         \blt[data generation]
         \begin{bullets}%
           \blt[scenario]The area of interest is a rectangle of
           $500\times 400$ m with 9 streets in each direction delimited by 8
           rows and columns of buildings of a certain height $\height$.
           \blt[fly grid]The flight height is between 50 and 150 m.
           \blt[slf]The SLF is such that the absorption inside the buildings is 3
           dB/m.  \blt[comms]The carrier frequency is 2.4 GHz, the bandwidth
           $\bandwidth=20$ MHz, the transmit power $\txpow=0.1$ Watt, and the
           noise power $\noisepow=-96$ dBm.
           \blt[users]A total of $\usernum$ GTs are deployed on the
           street uniformly at random.
         \end{bullets}%
         \blt[tested algorithms]
         \begin{bullets}%
           \blt[list]  The proposed algorithm is compared with the algorithm by
           \begin{bullets}%
             \blt[huang]Huang et al.~\cite{huang2020sparse},
             \blt[galkin]the
             K-means algorithm by Galkin et al.~\cite{galkin2016deployment},
             \blt[lyu]the spiral-based algorithm by Lyu et
             al.~\cite{lyu2017mounted},
             \blt[hammouti]and the iterative
             algorithm by Hammouti et al.~\cite{hammouti2019mechanism} for
             unlimited backhaul. 
           \end{bullets}%
           \blt[implementation] The implementation of the algorithm
           in~\cite{huang2020sparse} was provided by the authors, whereas the
           rest were implemented by us. The algorithm
           in~\cite{huang2020sparse} is only used in one experiment since its
           computational complexity of $\mathcal{O}(\usernum^6)$ makes it
           only suitable for a relatively low $\usernum$. The positions
           returned by these algorithms are projected onto the grid $\grid$
           of allowed flying positions.
         \end{bullets}%
         
         \blt[performance metrics]The adopted performance metric is the
         minimum number of ABSs required to guarantee  a rate $\ratemin$ to all
         GTs. This metric is averaged using  Monte Carlo  across
         realizations of the user locations.
         \blt[implementation]For the algorithms
         in~\cite{huang2020sparse} and~\cite{lyu2017mounted}, which
         are based on a maximum radius, the latter is gradually
         decreased starting from its value corresponding to free space
         propagation until all GTs receive the minimum rate. For the
         algorithms in~\cite{hammouti2019mechanism}
         and~\cite{galkin2016deployment}, the number of centroids is
         gradually increased starting from 1 until the aforementioned
         rate condition is met. See the repository (link on the first
         page) for more details along with the code of all
         experiments.

       \end{bullets}%

       \cmt{description of the experiments}
       \begin{bullets}%
         \blt[num users]Fig.~\ref{fig:num_users} depicts the minimum number
         of ABSs required to guarantee a rate of $\ratemin=5$ Mb/s for all
         GTs. The proposed algorithm is seen to yield placements that
         require fewer ABSs than all competing algorithms. This can be
         ascribed to the fact that it is aware of the channel and of in
         which regions it is allowed to fly.
         \blt[height]To investigate further the impact of the former
         effect, Fig.~\ref{fig:height} studies the influence of 
         shadowing. For a building height $\height=0$, propagation
         occurs in free space, which leads to all algorithms
         performing similarly. The slightly worse performance of the
         algorithm by Lyu et al. is mainly caused by the flight grid
         discretization. As $\height$ increases, the channel gradually
         differs more and more from free-space propagation and the
         competing algorithms suffer a performance degradation.
         
         \blt[rate]Finally, Fig.~\ref{fig:min_rate} investigates the influence of $\ratemin$. It
         is seen that the sensitivity of the proposed algorithm is much smaller
         than the one of its competitors, for which the performance metric
         increases considerably as $\ratemin$ increases.
         \begin{figure}[!t]
           \centering
           \includefig{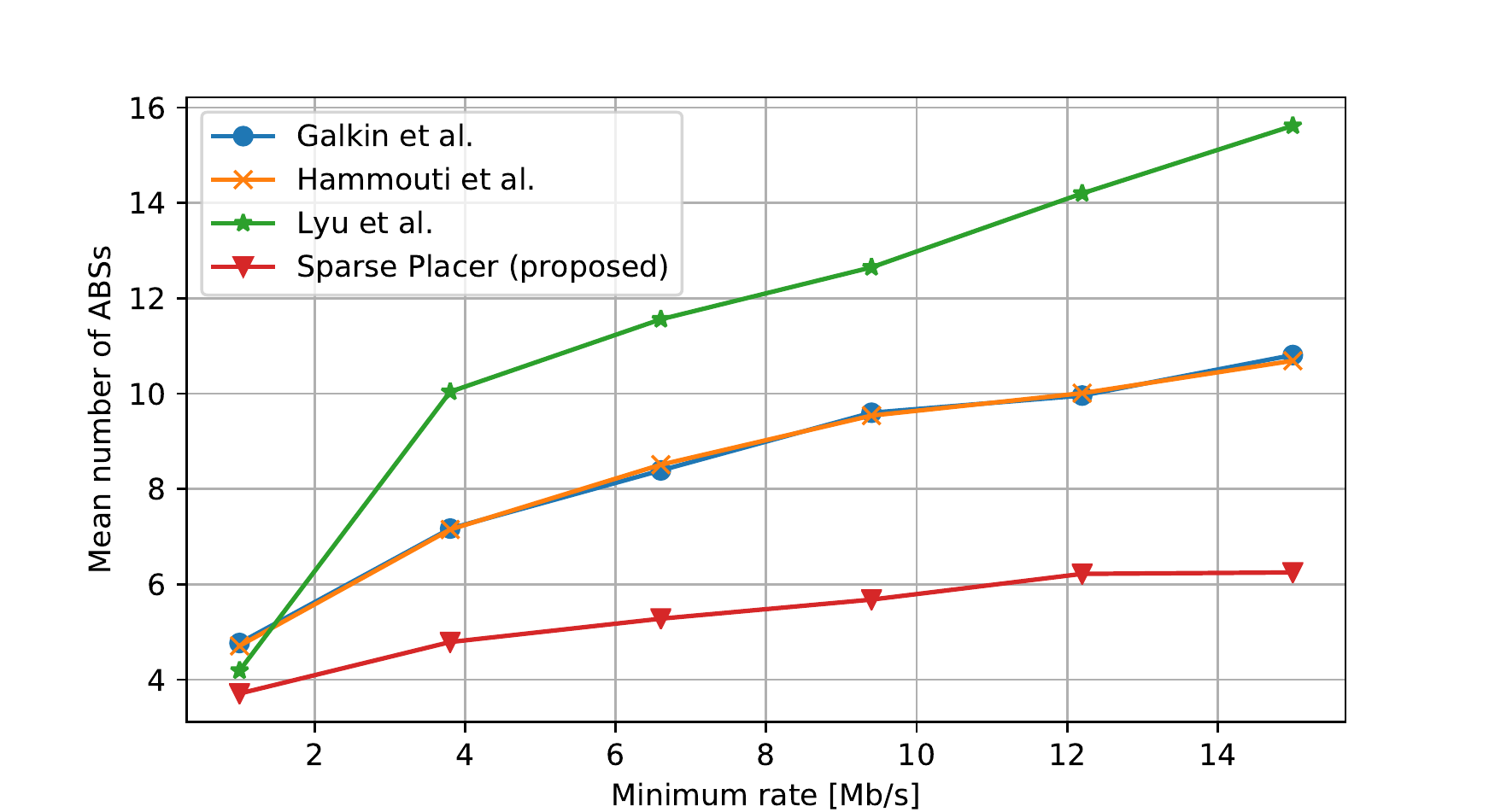}
           \caption{Mean minimum number of ABSs required to provide a minimum
             rate of $\ratemin$ ($\height=53$ m, $48\times 40\times 5$ SLF
             grid, $9\times 9\times 5$ fly grid).}
           \label{fig:min_rate}
         \end{figure}

       \end{bullets}%

       \section{Related Work}
       \label{sec:related}

       The most related works
       are~\cite{qiu2020reinforcement,chen2021relay,huang2020sparse}. 
       \begin{bullets}%
         \blt[qiu2020reinforcement]In~\cite{qiu2020reinforcement}, a
         terrain map or 3D model of the environment is used to predict
         the channel. 
         \begin{bullets}%
           \blt[limitations]
           \begin{bullets}%
             \blt[terrain map]Unfortunately, such models are seldom
             available and, furthermore, their resolution is typically very low
             relative to typical wavelengths, which indicates that the
             resulting accuracy may be insufficient for placement
             purposes. 
             \blt[DRL approach]Besides, a reinforcement learning approach
             is used rather than a convex optimization approach as in the
             present paper. 
             \begin{bullets}%
               \blt[] The algorithm needs to be retrained in every new
               environment or
               \blt[] if the number of UAVs changes.
               \blt[] Besides, this approach is not flexible enough to
               accommodate additional constraints, for example that a
               human user must take control of one of the UAVs.
             \end{bullets}%
           \end{bullets}%
         \end{bullets}%
         
         \blt[chen2021relay]The approach in~\cite{chen2021relay} relies
         on average local descriptors of the channel in terms of a map
         that provides the path loss exponent of each region in the
         deployment scenario.
         \begin{bullets}%
           \blt[1 uav]However, it just applies for $\absnum=\usernum=1$. 
           \blt[]\journal{
             a path loss exponent different from 2 is the result of
             fitting other phenomena besides free-space path loss, such as
             shadowing and fading. Thus, this model again captures propagation
             in average scenarios. }
         \end{bullets}%

         \blt[huang2020sparse]Finally, \cite{huang2020sparse} also adopts
         a convex optimization approach based on promoting sparsity, but
         the formulation is entirely different as it is not based on a
         discretization.
         \begin{bullets}%
           \blt[]Its complexity is $\mathcal{O}(\usernum^6)$,
           which restricts its applicability to scenarios with a low
           number of GTs.
           \blt[]Besides, it cannot accommodate general
           flight constraints since convexity would be lost in that case.
         \end{bullets}%

       \end{bullets}%

       \section{Conclusions}
       \label{sec:conclusions}

       This paper proposes a new approach to ABS placement where, instead of
       relying on average characterizations of the channel, a radio map of
       the specific deployment scenario is constructed and used to determine
       the set of optimal ABS locations in terms of a convex objective that
       approximately minimizes the number of ABSs to guarantee a minimum rate
       to all GTs. Unlike most approaches, the proposed algorithm has a low
       complexity and can accommodate flight constraints such as no-fly zones
       or airspace occupied by buildings. The intuitive soundness of the
       scheme is empirically corroborated using an open source simulator
       developed in this work.

\clearpage
\printmybibliography

\clearpage
\section{Supplementary Material}

\renewcommand{\journal}[1]{}

\subsection{Notation}
% If $\bm a$ is a column vector, then $a[m]$ denotes the $m$-th entry of
% $\bm a$.
$\rfield_{++}$ is the set of positive real numbers.
If $\bm a$ and $\bm b$ are vectors of the same dimension,
then $\bm a\odot \bm b$ is the entrywise product of $\bm a$ and
$\bm b$, whereas  $\bm a\div \bm b$ is the entrywise quotient of $\bm a$ and
$\bm b$.

\subsection{An Algorithm for Air-to-ground Radio Tomography}

\newcommand{\incs}{\hc{b}_\text{inc}}
\newcommand{\incvec}{\hc{\bm b}_\text{inc}}

\newcommand{\vis}{\hc{i}} % voxel index scalar
\newcommand{\vivec}{\hc{\bm i}} % voxel index vector
\newcommand{\vicurrents}{\hc{i}_\text{current}} 
\newcommand{\vicurrentvec}{\hc{\bm i}_\text{current}}
\newcommand{\viset}{\hc{\mathcal{I}}} % set of indices

\newcommand{\slftens}{\hc{L}}
\newcommand{\slften}{\hc{\bm \slftens}}

\newcommand{\deltaloc}{\hc{\bm \Delta}_{\loc}}
\newcommand{\deltagridvec}{\hc{\bm \delta}_{\slfgrid}}
\newcommand{\deltagrids}{\hc{ \delta}_{\slfgrid}}
\newcommand{\tcands}{\hc{ t}_\text{cand}}
\newcommand{\tcandvec}{\hc{\bm t}_\text{cand}}
\newcommand{\inext}{\hc{i}_\text{next}}
\newcommand{\tnext}{\hc{t}_\text{next}}
\newcommand{\integral}{\hc{{I}}}

  \begin{bullets}%
    \journal{
    \blt[std. approximation] Since no closed-form expression is
    generally available for $\slf$, the integral in \eqref{eq:tomoint}
    needs to be approximated numerically.
    \begin{bullets}%
      \blt[innner prod]  By far, the approximation
    used in the vast majority of works in radio tomography involves
    creating a regular  grid with points 
    $\{\slfgridpt_1,\ldots,\slfgridpt_\slfgridptnum\}$ and approximating
    \begin{align}
      \label{eq:tomointapprox}
      \shad(\loc_1,\loc_2)
      \approx \sum_\slfgridptind \weightfun(\loc_1, \loc_2, \slfgridpt_\slfgridptind)\slf(\slfgridpt_\slfgridptind)
    \end{align}
    \blt[weight fun] where the weight function
    $\weightfun(\loc_1, \loc_2, \slfgridpt)$ aims at assigning a non-zero
    weight only to those grid points lying close to the line segment
    between $\loc_1$ and $\loc_2$. This function typically takes a
    non-zero value only when $\slfgridpt$ lies inside an ellipse with
    foci at $\loc_1$ and $\loc_2$. A common example is
    \begin{align}
      \weightfun(\loc_1, \loc_2, \slfgridpt)=
      \begin{cases}
        0&\text{if } \|\loc_1-\slfgridpt\|_2 + \|\slfgridpt-\loc_2\|_2 \\
        &\quad >\|\loc_1 - \loc_2\| + \wfwidth/2\\
        1/{
      \| \loc_1 -\loc_2\|_2^{1/2}
} &\text{otherwise},
      \end{cases}
    \end{align}
    where $\wfwidth$ parameterizes the width of the ellipse and is
    typically set to the wavelength of the transmission since, in that
    way, the support of $\weightfun(\loc_1, \loc_2, \slfgridpt)$ is the
    first Fresnel ellipsoid. Other examples of can be found
    e.g. in~\cite{hamilton2014modeling,romero2018blind\journal{,estevezgutierrez2021hybrid}}.

    \blt[Limitations]The aforementioned approximation suffers from
    several limitations that render it impractical for the application
    at hand.
    \begin{bullets}%
      \blt[discontinuous]First, since
      $\weightfun(\loc_1, \loc_2, \slfgridpt)$ is discontinuous, the
      resulting approximation of $\shad(\loc_1,\loc_2)$ is also
      discontinuous and may lead to unpredictable behavior. For
      example, it may well happen that
      $\sum_\slfgridptind \weightfun(\loc_1, \loc_2,
      \slfgridpt_\slfgridptind)\slf(\slfgridpt_\slfgridptind)$ even when
      $\loc_1 \neq \loc_2$ and
      $\slf(\slfgridpt_\slfgridptind)~\forall \slfgridptind$. This happens when
      $\loc_1$ and $ \loc_2$ are such that no grid point falls inside
      the elliptical support of $\weightfun(\loc_1, \loc_2,
      \slfgridpt)$. \acom{figure?}
      \blt[large lambda] Thus, this approximation is more suitable to
      communications at low carrier frequencies, since in that case
      the wavelength is large and a large number of grid points are
      likely to be contained in the ellipse, thus mitigating these
      effects.       
      \blt[number of grid pts] But one could alternatively think that
      if the wavelength is short, then it suffices to create a
      sufficiently dense grid so that the distance between grid points
      is small, as compared to the wavelength. However, it is easy to
      see that this would entail a prohibitively high number of grid
      points. For example, if one considers a communication at 3 GHz,
      the wavelength is around 10 cm and one may want to set the grid
      points, say, 5 cm apart. If the region of interest is 1 km
      $\times$ 1 km $\times$ 100 m, then the total number of
      grid points is $10^{11}$, which is obviously impractical.
      \blt[complexity] Another critical limitation is computational
      complexity. Observe that \eqref{eq:tomointapprox} requires
      evaluating $\weightfun(\loc_1, \loc_2, \slfgridpt_\slfgridptind)$ and
      the product
      $ \weightfun(\loc_1, \loc_2,
      \slfgridpt_\slfgridptind)\slf(\slfgridpt_\slfgridptind)$ for each grid
      point. Thus, if the grid is
      $\slfgridside \times \slfgridside \times \slfgridside$, the complexity of
      approximating $\shad(\loc_1,\loc_2)$ is
      $\mathcal{O}(\slfgridside^3)$. This is computationally problematic
      for the application at hand since $\shad(\loc_1,\loc_2)$ is must
      be approximated at a large number of locations to solve
      \eqref{eq:problemf}.

    \end{bullets}%

  \end{bullets}%
} \blt[proposed]As indicated in Sec.~\ref{sec:radiomaps}, the usual
approximation to \eqref{eq:tomoint} using a weight function is not
suitable to construct an air-to-ground radio map for ABS placement.
  \begin{bullets}%
    \blt[overview]Instead, this work proposes adopting a different
    approximation to the integral in \eqref{eq:tomoint}. The
    technique, commonly used in other disciplines (see references in
    \cite{mitchell1990comparison}) and hinted in a different context
    in \cite{kanso2009compressed}, involves splitting the 3D space in
    voxels centered at the grid points
    $\slfgrid\define\{\slfgridpt_1,\ldots,\slfgridpt_\slfgridptnum\}$
    and approximating $\slf$ by a function that takes the value
    $\slf(\slfgridpt_\slfgridptind)$ at all points of the
    $\slfgridptind$-th voxel. The resulting piecewise constant
    approximation of $\slf$ can be integrated by determining the
    positions of the crossings between the voxel boundaries and the
    line segment between $\loc_1$ and $\loc_2$; see
    Fig.~\ref{fig:tomography}.

    \blt[algo explanation]Algorithm~\ref{algo:tomo}, which  can be classified as a
    parametric, floating point, and zeroth-order algorithm
    \cite[Sec. I-B-1]{mitchell1990comparison}, is our implementation 
    of the aforementioned approximation, yet  others are possible.
      \begin{bullets}%
        \blt[parameterization] The idea is to parameterize the line
        segment between $\loc_1$ and $\loc_2$ as
        $\loc(t)=\loc_1 + t(\loc_2-\loc_1)$, where $t\in[0,1]$, and
        identify the values $t_1,t_2,\ldots,t_T$ for which the
        boundary between two adjacent voxels is crossed.  \blt[sum]
        Since 
        $\|\loc(t_i)-\loc(t_{i-1})\| =
        (t_i-t_{i-1})\|\loc_2-\loc_1\|$ whenever  $t_i>t_{i-1}$, the approximation is then
    \begin{align}
      \shad(\loc_1,\loc_2)&\approx
                            \frac{\sum_{i=2}^T(t_i-t_{i-1})\|\loc_2-\loc_1\|
                            \slf(\slfgridpt_{\slfgridptind_i})}{\|\loc_2-\loc_1\|^{1/2}}\\
      &= \|\loc_2-\loc_1\|^{1/2}\sum_{i=2}^T(t_i-t_{i-1})\slf(\slfgridpt_{\slfgridptind_i}),
    \end{align}
    where ${\slfgridptind_i}$ is the index of the $i$-th voxel crossed by
    the segment.
    \blt[crossings]Since $\slfgrid$ is a 3D grid, each point in
    $\{\slfgridpt_1,\ldots,\slfgridpt_\slfgridptnum\}$ can also be
    indexed by a vector $\vivec$ of 3 indices that lies in the set
    $\viset\define\{1,\ldots,\slfgridsidex\}\times\{1,\ldots,\slfgridsidey\}\times\{1,\ldots,\slfgridsidez\}$. The
    values of the SLF can also be collected in a tensor
    $\slften\in\rfield^{\slfgridsidex\times \slfgridsidey\times
      \slfgridsidez}$, whose entry $\slftens[\vivec]$ is the value of
    $\slf$ at the $\vivec$-th grid point.  If
    $\deltagridvec\in \rfield_{++}^3$ denotes a vector whose $j$-th
    entry $\deltagrids\entnot{j}$ represents the spacing between grid
    points along the $j$-th axis, the coordinates of the $\vivec$-th
    grid point are clearly $\vivec\odot\deltagridvec$, where $\odot$
    denotes entrywise product. Similarly, the boundaries between
    adjacent voxels along the $j$-th axis occur at values of the
    $j$-th coordinate given by $\deltagrids\entnot{j}(\vis \pm 1/2)$,
    where $\vis$ is an integer.  It is then clear that
    steps~\ref{step:crossings}-\ref{step:nextcrossing} in
    Algorithm~\ref{algo:tomo} simply find the next value of $t$ for
    which the segment crosses a voxel boundary along one of the axes
    by solving the equation
    \begin{align}
      \locs_1\entnot{j} + t  (\locs_2\entnot{j}-\locs_1\entnot{j})=
      \deltagrids\entnot{j}(\vicurrents\entnot{j} \pm {1}/{2}) 
    \end{align}
    for $t$ along each axis $j$ and taking the minimum across
    axes. The $\pm$ becomes a plus sign for the $j$-th axis if the segment
    is increasing along this axis and a minus sign otherwise.

    \blt[alternative]An alternative implementation of the same
    integral approximation with smaller computational complexity but greater
    memory complexity could be obtained by creating 3 lists
    corresponding to the values of $t$ for which the line segment
    between $\loc_1$ and $\loc_2$ intersects each axis and then
    merging those lists into a list with non-decreasing values of $t$.
    %  it is clear that the values of $t$ for which the segment
    % crosses the boundaries of the voxel with indices
    % $\vicurrentvec\in \viset$ along the $j$-th axis satisfy

  \end{bullets}%

  \blt[algo benefits]Algorithm~\ref{algo:tomo} solves the
  limitations of the conventional approximation outlined in Sec.~\ref{sec:radiomaps}. 
    \begin{bullets}%
      \blt[continuous]First, Algorithm~\ref{algo:tomo} yields an
      approximation of $\shad(\loc_1,\loc_2)$ that is a continuous
      function of $\loc_1$ and $\loc_2$ since the line integral of a
      piecewise constant function is continuous. \journal{Besides, the
        issue of the approximation becoming zero when the elliptical
        support of the weight function in \eqref{eq:tomointapprox}
        misses all grid points now disappears.}\conference{Besides,
        the issue of the approximation becoming zero when the ellipses
        in the right side of Fig.~\ref{fig:tomography} miss all grid
        points disappears.}
      \blt[num. voxels]For this reason, the voxels can now be kept
      large regardless of the wavelength and, therefore, the total
      number of voxels can be kept low enough to be handled given the
      available computational resources.
      \blt[complexity]Finally, as indicated in
      Sec.~\ref{sec:radiomaps}, the computational complexity of
      Algorithm~\ref{algo:tomo} is much smaller than the one of the
      conventional approximation. Specifically, one can observe in
      Algorithm~\ref{algo:tomo} that a constant number of products and
      additions are required for each crossing. The total number of
      crossings is at most
      $\slfgridsidez+\slfgridsidey+\slfgridsidez$, which means that,
      if $\slfgridsidex=\slfgridsidey=\slfgridsidez=\slfgridside$,
      then the total complexity of Algorithm~\ref{algo:tomo} is
      $\mathcal{O}(\slfgridside)$, whereas the complexity of the
      standard approximation is $\mathcal{O}(\slfgridside^3)$.

    \end{bullets}%
  
  \end{bullets}%

  % \blt[tomography] \ra Bresenham's line algorithm \ra just determines through which pixels the line goes
\end{bullets}%

\begin{algorithm}[t!]               
\caption{Tomographic Integral Approximation}
\label{algo:tomo}    
\begin{minipage}{200cm}
\begin{algorithmic}[1]
  \STATE Input: $\loc_1$, $\loc_2$, grid spacing vector
  $\deltagridvec\in\rfield^3$, \\\quad\quad ~~SLF tensor
  $\slften\in\rfield^{\slfgridsidex\times \slfgridsidey\times \slfgridsidez}$. 
\STATE Initialize $\deltaloc=\loc_2-\loc_1$,  $\incvec =
\sign(\deltaloc)$, $\integral=0$
\STATE Set zero entries of $\deltaloc$ to 1 \# To avoid dividing by 0
\STATE Set $\vicurrentvec = \text{round}( \loc_1\div\deltagridvec)$ \#
Index of current voxel
\WHILE{$t<1$}

\STATE Set $\tcandvec =( \deltagridvec \odot (\vicurrentvec + \incvec/2) -
\loc_1)\div \deltaloc$
\label{step:crossings}

\STATE Set $\inext = \argmin_i
\tcands\entnot{i}~\st~\incs\entnot{i}\neq 0$

\STATE Set $\tnext = \tcands\entnot{\inext}$
\label{step:nextcrossing}

\STATE Set $\integral =\integral + (t-\tnext) \slftens\entnot{\vicurrentvec}$

\STATE Set $t = \tnext$

\STATE Set $\vicurrentvec\entnot{\inext} = \vicurrentvec\entnot{\inext} + \incs\entnot{\inext}$

\ENDWHILE
\STATE \textbf{return} $\|\loc_2-\loc_1\|^{1/2}\integral$
% \IF{convergence\_criterion() == True}
% \STATE return $\ratemat\itnot{\itind+1}$
% \ENDIF
\end{algorithmic}
\end{minipage}
\end{algorithm}

\subsection{Interior-Point Solver}
\label{sec:ip}
This section details how an interior-point solver can be used to solve
a relaxed version of \eqref{eq:probactisum}.

\begin{bullets}%
  \blt[relax constr.]Indeed, Problem \eqref{eq:probactisum} is
  non-convex due to the constraint
  $\activec\in \{ 0, 1\}^{ \gridptnum}$. As pointed out in
  Sec.~\ref{sec:placement}, a brute-force approach is not viable since
  $\gridptnum$ will typically be large in real applications. Instead,
  it is more convenient to adopt a convex approximation by relaxing
  this constraint. This yields
\begin{subequations}
  \begin{align}
    \minimize_{\activec\in[ 0, 1]^{ \gridptnum}}\quad &  \|\activec\|_0\\
    \st \quad & \capmat \activec \geq  \ratemin \bm 1 
  \end{align}
\end{subequations}
where $\activec\define[\acti_1,\ldots,\acti_\gridptnum]\transpose$ and
$\capmat\define[\capvec_1,\ldots,\capvec_\gridptnum]$.
\blt[relax obj.]Although it can be easily seen that this relaxation
does not entail loss of optimality, the objective now is non-convex.
As usual, this zero norm can be replaced with an $\ell_1$-norm to
yield a convex problem:
\begin{subequations}
  \begin{align}
    \minimize_{\activec\in[ 0, 1]^{ \gridptnum}}\quad &  \|\activec\|_1\\
    \st \quad & \capmat \activec \geq  \ratemin \bm 1.
  \end{align}
\end{subequations}
\blt[reweighting]It is well-known though that the sparsity of the
solutions can be increased by means of
reweighting~\cite{candes2008reweighted}. To this end,
$\|\activec\|_1=\sum_\gridptind |\acti_\gridptind|=\sum_\gridptind
\acti_\gridptind$ can be replaced with
$\sum_\gridptind \sparsweight_\gridptind\acti_\gridptind$, where 
$\sparsweight_\gridptind\geq 0$ are properly selected weights:
\begin{subequations}
  \begin{align}
    \minimize_{\activec\in[ 0, 1]^{ \gridptnum}}\quad &\sparsweightvec\transpose\activec\\
    \st \quad & \capmat \activec \geq  \ratemin \bm 1 %\\& \bm 0\leq \activec \leq \bm 1.
  \end{align}
\end{subequations}
The
standard approach is to iteratively set
$\sparsweight_\gridptind=1/(\epsilon + \actiold_\gridptind)$, where $\epsilon$ is a small constant and 
$\{\actiold_\gridptind\}_\gridptind$ are obtained by first solving the
problem with a previous set of weights, the initial set being such that $\sparsweight_\gridptind=1~\forall\gridptind$.

\blt[interior pt]To apply an interior point solver, the inequality
constraints can be replaced with equality constraints by introducing
the vector of slack variables $\slackvec$:
\begin{subequations}
  \begin{align}
    \minimize_{\activec\in[ 0, 1]^{ \gridptnum}}\quad &\sparsweightvec\transpose\activec\\
    \st \quad & \capmat \activec =  \ratemin \bm 1 + \slackvec\\ %\\& \bm 0\leq \activec \leq \bm 1.
    & \slackvec \geq \bm 0
  \end{align}
\end{subequations}
Given that this problem has $\gridptnum + \usernum$ variables plus
$\usernum$ Lagrange multipliers associated with the equality
constraints, each inner step of the interior point solver involves
solving a system of $\gridptnum + 2\usernum$ linear equations with
$\gridptnum + 2\usernum$ variables, which has a complexity
$\mathcal{O}((\gridptnum + 2\usernum)^3)$.

\end{bullets}%

%\ref{algo:placement}

\subsection{Derivation of Algorithm \ref{algo:placement}}

The goal in this section is to derive a solver for
\eqref{eq:ubeqproblem} using the framework of
ADMM~\cite{boyd2011distributed}.
\begin{bullets}%
%   \blt[Preparing the problem]
%   \begin{bullets}%
%     \blt[Objective] The second step is to replace the objective
%     function with a convex surrogate. Clearly, the objective in
%     \eqref{eq:ubproblem} equals the zero norm of the vector
%     $\normvec\define[\|\ratevec_1\|_\infty, \ldots,
%     \|\ratevec_\gridptnum\|_\infty]\transpose$. As in
%     Sec.~\ref{sec:ip}, one can follow the conventional approach of
%     replacing this zero norm with an $\ell_1$ norm
%     $\|\normvec\|_1\define \sum_\gridptind |
%     \norm\vecent{\gridptind}|= \sum_\gridptind
%     \|\ratevec_\gridptind\|_\infty $. More precisely, the approach
%     here is to replace it with a weighted $\ell_1$-norm loss of the
%     form
%     $\sum_\gridptind \sparsweight_\gridptind
%     \|\ratevec_\gridptind\|_\infty $, where the non-negative weights are set as in
%     \acom{cite}.
% \end{bullets}%

  \blt[slack variables]To facilitate this task, it is convenient at
  this point to replace the objective with a linear function
  $\sum_\gridptind \sparsweight_\gridptind \slack_\gridptind$, where
  $\{\slack_\gridptind\}_\gridptind$ are slack variables. Problem
  \eqref{eq:ubeqproblem} can then be expressed as
\begin{subequations}
  \label{eq:ubslack}
  \begin{align}
  \minimize_{\ratemat\in \rfield^{\usernum \times \gridptnum}}\quad & 
    \sparsweightvec\transpose \slackvec\\
  \st \quad & \ratemat \bm 1 = \ratemin \bm 1 \\
                                                                    & \bm 0 \leq \ratemat \leq \capmat\\
                                                                    & \ratevec_\gridptind\leq \slack_\gridptind \bm 1,~\gridptind = 1,\ldots,\gridptnum,
\end{align}
\end{subequations}
where
$\sparsweightvec\define[\sparsweight_1,\ldots,\sparsweight_\gridptnum]\transpose$
and $\slackvec\define [\slack_1,\ldots,\slack_\gridptnum]\transpose$.

\blt[formulation in ADMM form]The next step is to
express \eqref{eq:ubslack} in form amenable to application of ADMM.
  \begin{bullets}%
    \blt[Theory]For the problem at hand, notation can be simplified by
    adopting the following special homogeneous form:
\begin{subequations}
  \label{eq:admmstd}
  \begin{align}
  \minimize_{\Xmat, \Zmat}\quad & \admmfunx(\Xmat)+\admmfunz(\Zmat)\\
    \st \quad & \Amat_1\Xmat \Amat_2+ \Bmat_1\Zmat\Bmat_2 = \bm 0.
  \end{align}
\end{subequations}

For this problem, the ADMM iteration from
\cite[Sec.~3.1.1]{boyd2011distributed} becomes
\begin{subequations}
  \label{eq:admmstdit}
  \begin{align}
    \label{eq:admmstditx}
    \Xmat\itnot{\itind+1} &= \argmin_\Xmat \admmfunx(\Xmat) + \frac{\admmstep}{2}\| \Amat_1\Xmat \Amat_2+ \Bmat_1\Zmat\itnot{\itind}\Bmat_2 + \Umat\itnot{\itind}\|_\frob^2\\
        \label{eq:admmstditz}
    \Zmat\itnot{\itind+1} &= \argmin_\Zmat \admmfunz(\Zmat) + \frac{\admmstep}{2}\| \Amat_1\Xmat\itnot{\itind+1} \Amat_2+ \Bmat_1\Zmat\Bmat_2 + \Umat\itnot{\itind}\|_\frob^2\\
    \label{eq:admmstditu}
    \Umat\itnot{\itind+1} &=     \Umat\itnot{\itind} + \Amat_1\Xmat\itnot{\itind+1} \Amat_2+ \Bmat_1\Zmat\itnot{\itind+1}\Bmat_2,
  \end{align}
\end{subequations}
where $\Umat\itnot{\itind}$ is a  matrix of scaled dual variables and
$\admmstep>0$ is the step-size parameter.

\blt[assignments] There are multiple possibilities to cast
\eqref{eq:ubslack} as \eqref{eq:admmstd} and each one leads to updates
of a different nature. Thus, several attempts are often required. As
seen later, the following assignments yield suitable updates for the
problem under consideration:
\begin{subequations}
  \label{eq:ubassignments}
  \begin{align}
    \Xmat &\rightarrow [\ratemat\transpose , \slackvec]\transpose\\
    \Zmat &\rightarrow \ratemat\\
    \admmfunx(\Xmat)&\rightarrow     \sparsweightvec\transpose \slackvec + \sum_\gridptind
            \indicatorinf[\ratevec_\gridptind\leq \slack_\gridptind \bm 1]\\
    \admmfunz(\Xmat)&\rightarrow \indicatorinf[\ratemat \bm 1 = \ratemin \bm 1] + \indicatorinf[\bm 0 \leq \ratemat \leq \capmat]\\
    \Amat_1 &\rightarrow [\bm I_\usernum, \bm 0], ~\Amat_2 \rightarrow \bm I_{\gridptnum }, ~\Bmat_1\rightarrow - \bm I_\usernum, ~\Bmat_2\rightarrow \bm I_\gridptnum.
  \end{align}
\end{subequations}
Here, $\indicatorinf[\cdot]$ is a function that takes the value 0 when
the condition inside brackets holds and $\infty$ otherwise. Note that
with this choice for the matrices in \eqref{eq:admmstd}, it follows
that $\Amat_1\Xmat \Amat_2+ \Bmat_1\Zmat\Bmat_2=\ratemat-\Zmat$ and, therefore, the
constraint in \eqref{eq:admmstd} imposes that $\ratemat=\Zmat$.

  \end{bullets}%

  \blt[X-step]\textbf{$\Xmat$-step.}
  \begin{bullets}%
    \blt[Special expression]
  To derive the $\Xmat$-update, observe that, with the above
  assignments, the problem in \eqref{eq:admmstditx} becomes
\begin{subequations}
  \begin{align}
    (\ratemat\itnot{\itind+1},\slackvec\itnot{\itind+1}) =& \arg\min_{\ratemat, \slackvec}~
                                          \sparsweightvec\transpose \slackvec + \sum_\gridptind
            \indicatorinf[\ratevec_\gridptind\leq \slack_\gridptind \bm 1]
                                          \nonumber\\ &+ \frac{\admmstep}{2}\|\ratemat - \Zmat\itnot{\itind
                                                           }+\Umat\itnot{\itind}\|_\frob^2
    \\
    =& \arg\min_{\ratemat, \slackvec}~
                                        \sum_\gridptind \big[ \sparsweight_\gridptind \slack_\gridptind + 
            \indicatorinf[\ratevec_\gridptind\leq \slack_\gridptind \bm 1]
                                          \nonumber\\&+ \frac{\admmstep}{2}\|\ratevec_\gridptind - \zvec_\gridptind\itnot{\itind
                                          }+\uvec_\gridptind\itnot{\itind}\|_2^2                                          \big],
  \end{align}
\end{subequations}
where $\zvec_\gridptind\itnot{\itind }$ and $\uvec_\gridptind\itnot{\itind}$
respectively denote the $\gridptind$-th column of $\Zmat\itnot{\itind}$
and $\Umat\itnot{\itind}$.
\blt[separation] This problem clearly separates into
$\gridptnum$ problems of the form
\begin{subequations}
  \label{eq:ratevecindividual}
  \begin{align}
    (\ratevec_\gridptind\itnot{\itind+1},\slack_\gridptind\itnot{\itind+1}) 
    = \arg\min_{\ratevec_\gridptind, \slack_\gridptind}~& \sparsweight_\gridptind \slack_\gridptind 
                                          + \frac{\admmstep}{2}\|\ratevec_\gridptind - \zvec_\gridptind\itnot{\itind
      }+\uvec_\gridptind\itnot{\itind}\|_2^2 \\
    \st\quad&       \ratevec_\gridptind\leq \slack_\gridptind \bm 1.
  \end{align}
\end{subequations}
\blt[nec. and suff. conds] If $\sparsweight_\gridptind=0$, the inequality constraint can
be removed and the optimum is attained when
$\ratevec_\gridptind\itnot{\itind+1}= \zvec_\gridptind\itnot{\itind
}-\uvec_\gridptind\itnot{\itind}$. Thus, it suffices to focus on the case
$\sparsweight_\gridptind>0$. In this case, we have the following:

\begin{myproposition}
  If $\sparsweight_\gridptind>0$, then
  $\ratevec_\gridptind\itnot{\itind+1}$ and
  $\slack_\gridptind\itnot{\itind+1} $ satisfy
  \begin{subequations}
    \label{eq:xkkt}
    \begin{align}
      \label{eq:xkktr}
    &\ratevec_\gridptind\itnot{\itind+1} = \min(\zvec_\gridptind\itnot{\itind
      }-\uvec_\gridptind\itnot{\itind}, \slack_\gridptind\itnot{\itind+1}\bm 1)\\
      \label{eq:xkks}
    &\bm 1\transpose \max(\zvec_\gridptind\itnot{\itind
    }-\uvec_\gridptind\itnot{\itind} -\slack_\gridptind\itnot{\itind+1}\bm 1, \bm 0) = \frac{\sparsweight_\gridptind}{\admmstep}, 
  \end{align}
\end{subequations}
where  $\min$ and $\max$ operate  entrywise.
\end{myproposition}  

\begin{IEEEproof}
  Since Problem \eqref{eq:ratevecindividual} is convex differentiable
  and Slater's conditions are satisfied, it follows that the
  Karush-Kuhn-Tucker (KKT) conditions are sufficient and
  necessary\journal{~\cite[Sec.~5.5.3]{boyd}}. To obtain these
  conditions, observe that the Lagrangian of
  \eqref{eq:ratevecindividual} is given by
  \begin{align}
    \lagrangian(\ratevec_\gridptind, \slack_\gridptind;\nuvec)
    = \sparsweight_\gridptind \slack_\gridptind 
                                          + \frac{\admmstep}{2}\|\ratevec_\gridptind - \zvec_\gridptind\itnot{\itind
      }+\uvec_\gridptind\itnot{\itind}\|_2^2 + \nuvec\transpose(   \ratevec_\gridptind- \slack_\gridptind \bm 1).
  \end{align}
The KKT conditions are, therefore, 
\begin{subequations}
      \label{eq:ubkkt}
      \begin{align}
            \label{eq:ubkktgrate}
    \nabla_{\ratevec_\gridptind}\lagrangian(\ratevec_\gridptind, \slack_\gridptind;\nuvec) &=
 \admmstep(\ratevec_\gridptind - \zvec_\gridptind\itnot{\itind
                                                                                           }+\uvec_\gridptind\itnot{\itind}) + \nuvec = \bm 0\\
                \label{eq:ubkktslack}
    \nabla_{\slack_\gridptind}\lagrangian(\ratevec_\gridptind, \slack_\gridptind;\nuvec) &=
                                                                                         \sparsweight_\gridptind  - \bm 1\transpose \nuvec=0\\
        \label{eq:ubkktprimal}
        \ratevec_\gridptind&\leq \slack_\gridptind \bm 1\\
        \label{eq:ubkktnu}
    \nuvec&\geq \bm 0,~~\nus\entnot{\userind}(\rate_\gridptind\entnot{\userind} - \slack_\gridptind) = 0~\forall \userind.   
  \end{align}
\end{subequations}

From \eqref{eq:ubkktgrate} and the inequality in \eqref{eq:ubkktnu}, it follows that
\begin{align}
  \label{eq:ubkktvnucomb}
   \nuvec = -\admmstep(\ratevec_\gridptind - \zvec_\gridptind\itnot{\itind
  }+\uvec_\gridptind\itnot{\itind})\geq \bm 0.
  \end{align}
  This implies that
  $ \ratevec_\gridptind\leq \zvec_\gridptind\itnot{\itind } -
  \uvec_\gridptind\itnot{\itind}$.  Combining this inequality with
  \eqref{eq:ubkktprimal} yields
  \begin{align}
    \label{eq:ubkktratebound}
   \ratevec_\gridptind\leq \min(\zvec_\gridptind\itnot{\itind } -
  \uvec_\gridptind\itnot{\itind}, \slack_\gridptind \bm 1).
  \end{align}
 On the other hand, from the
  equality in \eqref{eq:ubkktvnucomb} and the inequality in
  \eqref{eq:ubkktnu}, one finds that
  \begin{align}
-\admmstep(\rate_\gridptind\entnot{\userind} - \zs_\gridptind\itnot{\itind
  }\entnot{\userind}+\us_\gridptind\itnot{\itind}\entnot{\userind})(\rate_\gridptind\entnot{\userind} - \slack_\gridptind) = 0~\forall \userind.
  \end{align}
This holds if and only if either $\rate_\gridptind\entnot{\userind} = \zs_\gridptind\itnot{\itind
  }\entnot{\userind}-\us_\gridptind\itnot{\itind}\entnot{\userind}$ or 
  $\rate_\gridptind\entnot{\userind} = \slack_\gridptind$. Therefore, it follows from \eqref{eq:ubkktratebound} that
    \begin{align}
    \label{eq:ubkktrateeq}
   \ratevec_\gridptind= \min(\zvec_\gridptind\itnot{\itind } -
  \uvec_\gridptind\itnot{\itind}, \slack_\gridptind \bm 1),
  \end{align}
  which establishes \eqref{eq:xkktr}. Finally, combine this expression with
  \eqref{eq:ubkktslack} and \eqref{eq:ubkktvnucomb} to arrive at
  \begin{subequations}
  \begin{align}
    \sparsweight_\gridptind  &=  -\admmstep \bm 1\transpose(\ratevec_\gridptind - \zvec_\gridptind\itnot{\itind
                               }+\uvec_\gridptind\itnot{\itind})\\
    &=  -\admmstep \bm 1\transpose(\min(\zvec_\gridptind\itnot{\itind } -
  \uvec_\gridptind\itnot{\itind}, \slack_\gridptind \bm 1) - \zvec_\gridptind\itnot{\itind
  }+\uvec_\gridptind\itnot{\itind})\\
    &=  -\admmstep \bm 1\transpose\min(\bm 0, \slack_\gridptind \bm 1 - \zvec_\gridptind\itnot{\itind
      }+\uvec_\gridptind\itnot{\itind}) \\
    \label{eq:ubslackeqp}
    &=  \admmstep \bm 1\transpose\max(\bm 0,  \zvec_\gridptind\itnot{\itind
  }-\uvec_\gridptind\itnot{\itind}-\slack_\gridptind \bm 1 ),
  \end{align}
\end{subequations}
thereby recovering \eqref{eq:xkks}. The proof is complete by noting
that \eqref{eq:ubkkt} holds if and only if \eqref{eq:ubkktrateeq} and
\eqref{eq:ubslackeqp} hold. 
  
\end{IEEEproof}

\blt[solving the conditions]
Observe that \eqref{eq:xkktr} can be used to obtain
$\ratevec_\gridptind\itnot{\itind+1}$ if
$\slack_\gridptind\itnot{\itind+1}$ is given, whereas \eqref{eq:xkks}
does not depend on $\ratevec_\gridptind\itnot{\itind+1}$. Therefore, a
solution to \eqref{eq:xkkt} can be found by first solving
\eqref{eq:xkks} for $\slack_\gridptind\itnot{\itind+1}$ and then
substituting the result into the right-hand side of \eqref{eq:xkktr}
to recover $\ratevec_\gridptind\itnot{\itind+1}$. To this end, we have
the following:

\begin{myproposition}
  \thlabel{prop:rootssg}
  Equation \eqref{eq:xkks} has a unique root. This root lies in the interval $[\slacklow\gii, \slackhigh\gii]$, where
  \begin{subequations}
    \begin{align}
      \label{eq:slacklow}
    \slacklow\gii &\define \min_\userind\left( \zs\gii\entnot{\userind} - \us\gii\entnot{\userind}\right) - \frac{\sparsweight_\gridptind}{\usernum \admmstep}\\
    \label{eq:slackhigh}
        \slackhigh\gii &\define \max_\userind\left( \zs\gii\entnot{\userind} - \us\gii\entnot{\userind}\right) - \frac{\sparsweight_\gridptind}{\usernum \admmstep}.
  \end{align}
\end{subequations}

\end{myproposition}

\begin{IEEEproof}
  Consider the function
  $\Ffun(\slack)\define \bm 1\transpose \max(\zs_\gridptind\itnot{\itind
  }-\us_\gridptind\itnot{\itind} -\slack\bm 1, \bm 0) = \sum_\userind
  \max(\zs_\gridptind\itnot{\itind
  }\entnot{\userind}-\us_\gridptind\itnot{\itind}\entnot{\userind}
  -\slack, \bm 0)$. Since $\Ffun$ is the sum of non-increasing piecewise
  linear functions, so is $\Ffun$. Since $\Ffun(\slack)\rightarrow \infty$
  as $\slack\rightarrow -\infty$ and $\Ffun(\slack)=0$ for a sufficiently
  large $\slack$, it follows that \eqref{eq:xkks} has at least one
  root. Uniqueness of the root follows readily by noting that 
  $\Ffun$ is strictly decreasing  whenever $\Ffun(\slack)>0$.

  It remains to be shown that
  $\Ffun(\slacklow\gii)\geq{\sparsweight_\gridptind}/({\usernum
    \admmstep})$ whereas
  $\Ffun(\slackhigh\gii)\leq {\sparsweight_\gridptind}/({\usernum
    \admmstep})$.
  \begin{bullets}%
    \blt[lower]
  For the first of these inequalities, observe that
  $\slacklow\gii \leq \zs\gii\entnot{\userind} -
  \us\gii\entnot{\userind} - {\sparsweight_\gridptind}/({\usernum
    \admmstep})$ for all $\userind$, which in turn implies that 
  $\zs\gii\entnot{\userind} -
  \us\gii\entnot{\userind} - \slacklow\gii \geq {\sparsweight_\gridptind}/({\usernum
    \admmstep})$. Thus, $\max(\zs\gii\entnot{\userind} -
  \us\gii\entnot{\userind} - \slacklow\gii,0)=\zs\gii\entnot{\userind} -
  \us\gii\entnot{\userind} - \slacklow\gii \geq {\sparsweight_\gridptind}/({\usernum
    \admmstep})$, which yields $\Ffun(\slacklow\gii)\geq \sum_\userind{\sparsweight_\gridptind}/({\usernum
    \admmstep}) = {\sparsweight_\gridptind}/{
    \admmstep}$.
  \blt[upper]For the second inequality, note similarly that   $\zs\gii\entnot{\userind} -
  \us\gii\entnot{\userind} - \slackhigh\gii \leq {\sparsweight_\gridptind}/({\usernum
    \admmstep})$ for all $\userind$. This means that $\Ffun(\slackhigh\gii)\leq \sum_\userind\max(
    {\sparsweight_\gridptind}/({\usernum
    \admmstep}),0) =  {\sparsweight_\gridptind}/{
    \admmstep} $.

  \end{bullets}%
  
\end{IEEEproof}

\end{bullets}%

\blt[Z-step]\textbf{$\Zmat$-step.}
  \begin{bullets}%
    \blt[Specialization]From \eqref{eq:admmstditz} and \eqref{eq:ubassignments}, it follows that
    \begin{subequations}
  \begin{align}
    \Zmat\itnot{\itind+1} = \argmin_\Zmat ~~&\bigg[ \indicatorinf[\Zmat \bm 1 = \ratemin \bm 1] + \indicatorinf[\bm 0 \leq \Zmat \leq \capmat] \\&+ \frac{\admmstep}{2}\| \ratemat\itnot{\itind+1} -\Zmat + \Umat\itnot{\itind}\|_\frob^2\bigg]\\
     = \argmin_\Zmat ~~&\sum_\userind\bigg[ \indicatorinf[\zrvec_\userind\transpose \bm 1 = \ratemin] + \indicatorinf[\bm 0 \leq \zrvec_\userind \leq \caprvec_\userind] \\&+ \frac{\admmstep}{2}\| \ratervec_\userind\itnot{\itind+1} -\zrvec_\userind + \urvec_\userind\itnot{\itind}\|_\frob^2\bigg],
  \end{align}
\end{subequations}
where $\zrvec_\userind$, $\caprvec_\userind$, and $\urvec_\userind\itnot{\itind}$ respectively denote the $\userind$-th column of $\Zmat\transpose$, $\capmat\transpose$, and $(\Umat\itnot{\itind})\transpose$. 
\blt[separation] Clearly, this separates into $\usernum$ problems of the form
\begin{subequations}
  \label{eq:ubzrvecs}
  \begin{align}
    \zrvec_\userind\itnot{\itind+1} = \argmin_{\zrvec_\userind} ~~&  \frac{1}{2}\| \ratervec_\userind\itnot{\itind+1} -\zrvec_\userind + \urvec_\userind\itnot{\itind}\|_\frob^2\\
    \st \quad &  \bm 1\transpose \zrvec_\userind = \ratemin,~~ \bm 0 \leq \zrvec_\userind \leq \caprvec_\userind.
  \end{align}
\end{subequations}

\blt[Nec. and suff. conds]
\begin{myproposition} If
  $\bm 1\transpose \caprvec_\userind< \ratemin$, then
  \eqref{eq:ubzrvecs} is infeasible.  If
  $\bm 1\transpose \caprvec_\userind\geq \ratemin$, the solution to
  \eqref{eq:ubzrvecs} is given by
  \begin{align}
    \label{eq:ubzrvecsol}
    \zrvec_\userind\itnot{\itind+1} =
\max(\bm 0,    \min(\caprvec_\userind, \ratervec_\userind\itnot{\itind+1} +
    \urvec_\userind\itnot{\itind} - \lambdas\bm 1)),
  \end{align}
  where $\lambdas$ satisfies
  \begin{align}
    \label{eq:ubzrveclambda}
\bm 1\transpose\max(\bm 0,    \min(\caprvec_\userind, \ratervec_\userind\itnot{\itind+1} +
    \urvec_\userind\itnot{\itind} - \lambdas\bm 1)) = \ratemin.
  \end{align}

\end{myproposition}

\begin{IEEEproof}
  The fact that $\bm 1\transpose \caprvec_\userind< \ratemin$ implies
  that \eqref{eq:ubzrvecs} is infeasible is trivial and, therefore,
  the rest of the proof focuses on the case where
  $\bm 1\transpose \caprvec_\userind\geq \ratemin$.
  
  As before, the KKT conditions are sufficient and necessary in this
  case. Noting that  the Lagrangian is given by
  \begin{align}
    &\lagrangian(\zrvec_\userind;\lambdas, \nuvec,\muvec ) =~
     \frac{1}{2}\| \ratervec_\userind\itnot{\itind+1} -\zrvec_\userind + \urvec_\userind\itnot{\itind}\|_\frob^2 \nonumber\\&\quad\quad+ \lambdas( \bm 1\transpose \zrvec_\userind - \ratemin) - \nuvec\transpose \zrvec_\userind + \muvec\transpose( \zrvec_\userind - \caprvec_\userind)
  \end{align}
  yields the KKT conditions
  \begin{subequations}
    \begin{align}
      \nonumber
      &\nabla_{\zrvec_\userind}    \lagrangian(\zrvec_\userind;\lambdas, \nuvec,\muvec ) =\\
      \label{eq:ubkktz}
    &\quad
      -( \ratervec_\userind\itnot{\itind+1} -\zrvec_\userind + \urvec_\userind\itnot{\itind}) + \lambdas \bm 1 - \nuvec + \muvec = \bm 0,\\
                  \label{eq:ubkktzlambda}
      &\bm 1\transpose \zrvec_\userind = \ratemin, \\
            \label{eq:ubkktznu}
      & \zrvec_\userind \geq \bm 0,~\nuvec \geq 0,~\nus\entnot{\gridptind}\zrs_\userind\entnot{\gridptind}=0~\forall \gridptind,\\
      \label{eq:ubkktmu}
      & \zrvec_\userind \leq \caprvec_\userind,~\muvec \geq 0,~\mus\entnot{\gridptind}(\zrs_\userind\entnot{\gridptind} -\caprs_\userind\entnot{\gridptind} )=0~\forall \gridptind.
  \end{align}
\end{subequations}
From \eqref{eq:ubkktz} and the second inequality in \eqref{eq:ubkktmu}, it follows that
\begin{align}
  \label{eq:ubmuvecin}
  \muvec =  \ratervec_\userind\itnot{\itind+1} -\zrvec_\userind + \urvec_\userind\itnot{\itind}
  - \lambdas \bm 1 + \nuvec \geq \bm 0,
  \end{align}
which in turn implies that 
\begin{align}
\zrvec_\userind \leq  \ratervec_\userind\itnot{\itind+1} + \urvec_\userind\itnot{\itind}
  - \lambdas \bm 1 + \nuvec .
\end{align}
Combining this expression with the first inequality in \eqref{eq:ubkktmu} yields
\begin{align}
\zrvec_\userind \leq \min(\caprvec_\userind,  \ratervec_\userind\itnot{\itind+1} + \urvec_\userind\itnot{\itind}
  - \lambdas \bm 1 + \nuvec ).
\end{align}
To show that this expression holds with equality, substitute \eqref{eq:ubmuvecin} into the equality of \eqref{eq:ubkktmu} to obtain
\begin{align}
  (\raters_\userind\itnot{\itind+1}\entnot{\gridptind} -\zrs_\userind\entnot{\gridptind} + \urs_\userind\itnot{\itind}\entnot{\gridptind}
  - \lambdas  + \nus\entnot{\gridptind}
  )(\zrs_\userind\entnot{\gridptind} -\caprs_\userind\entnot{\gridptind} )=0, 
  \end{align}
  which implies that either $ \zrs_\userind\entnot{\gridptind} = \raters_\userind\itnot{\itind+1}\entnot{\gridptind}  + \urs_\userind\itnot{\itind}\entnot{\gridptind}
  - \lambdas  + \nus\entnot{\gridptind}$ or $\zrs_\userind\entnot{\gridptind} =\caprs_\userind\entnot{\gridptind}$. Therefore,
  \begin{align}
    \label{eq:ubzrvec}
\zrvec_\userind = \min(\caprvec_\userind,  \ratervec_\userind\itnot{\itind+1} + \urvec_\userind\itnot{\itind}
  - \lambdas \bm 1 + \nuvec ).
  \end{align}
  To obtain an expression for $\zrvec_\userind$ that does not depend on 
 $\nuvec$, one may consider three cases for each $\gridptind$:
  \begin{itemize}
    \item C1:
    $\raters_\userind\itnot{\itind+1}\entnot{\gridptind} +
    \urs_\userind\itnot{\itind}\entnot{\gridptind} - \lambdas< 0$. In
    this case, if $\nus\entnot{\gridptind}=0$, expression
    \eqref{eq:ubzrvec} would imply that
    $\zrs_\userind\entnot{\gridptind}<0$, which would violate the
    first inequality in \eqref{eq:ubkktznu}. Therefore,
    $\nus\entnot{\gridptind}>0$ and, due to the equality in
    \eqref{eq:ubkktznu}, $\zrs_\userind\entnot{\gridptind}=0$. If
    $\caprs_\userind\entnot{\gridptind}>0$, it is then clear from
    \eqref{eq:ubzrvec} that
    $ \nus\entnot{\gridptind} =
    -(\raters_\userind\itnot{\itind+1}\entnot{\gridptind} +
    \urs_\userind\itnot{\itind}\entnot{\gridptind} - \lambdas)$.  If
    $\caprs_\userind\entnot{\gridptind}=0$, then greater values of
    $\nus\entnot{\gridptind}$ will also satisfy the KKT conditions but
    this is not relevant since in this case the only feasible
    $\zrs_\userind\entnot{\gridptind}$ is
    $\zrs_\userind\entnot{\gridptind}=0$.
    
    \item C2:
    $\raters_\userind\itnot{\itind+1}\entnot{\gridptind} +
    \urs_\userind\itnot{\itind}\entnot{\gridptind} - \lambdas = 0$. In
    this case, \eqref{eq:ubzrvec} becomes
    $\zrs_\userind\entnot{\gridptind} =
    \min(\caprs_\userind\entnot{\gridptind},
    \nus\entnot{\gridptind})$. Due to the equality in
    \eqref{eq:ubkktznu}, it then follows that either
    $\caprs_\userind\entnot{\gridptind}=0$ and
    $\nus\entnot{\gridptind}\geq 0$, or
    $\zrs_\userind\entnot{\gridptind} = \nus\entnot{\gridptind}=0$.

    \item C3:
    $\raters_\userind\itnot{\itind+1}\entnot{\gridptind} +
    \urs_\userind\itnot{\itind}\entnot{\gridptind} - \lambdas>0$. If
    $\caprs_\userind\entnot{\gridptind}=0$, then necessarily
    $\zrs_\userind\entnot{\gridptind} = 0$ and any
    $\nus\entnot{\gridptind}\geq 0$ satisfies the KKT conditions. On
    the other hand, if $\caprs_\userind\entnot{\gridptind}>0$, then it
    is clear that $\zrs_\userind\entnot{\gridptind} > 0$ and, due to
    the equality in \eqref{eq:ubkktznu}, one has that
    $\nus\entnot{\gridptind} =0$, which in turn implies that
    $\zrs_\userind\entnot{\gridptind} =
    \min(\caprs_\userind\entnot{\gridptind}, \raters_\userind\itnot{\itind+1}\entnot{\gridptind} +
    \urs_\userind\itnot{\itind}\entnot{\gridptind} - \lambdas)$.
  \end{itemize}
  Combining C1-C3 yields
  \begin{align}
    \label{eq:ubcombzrs}
    \zrs_\userind\entnot{\gridptind} =
\max(0,    \min(\caprs_\userind\entnot{\gridptind}, \raters_\userind\itnot{\itind+1}\entnot{\gridptind} +
    \urs_\userind\itnot{\itind}\entnot{\gridptind} - \lambdas)), 
    \end{align}
    which is just the scalar version of
    \eqref{eq:ubzrvecsol}. Finally, to obtain $\lambdas$, one may
    substitute \eqref{eq:ubcombzrs} into \eqref{eq:ubkktzlambda},
    which produces \eqref{eq:ubzrveclambda}. 

\end{IEEEproof}

\blt[solving the conditions] Thus, as in the $\Xmat$-step, one needs
to solve the scalar equation \eqref{eq:ubzrveclambda}. The following
result is the counterpart of \thref{prop:rootssg} for the
$\Zmat$-step.

\begin{myproposition}
  \thlabel{prop:ubrootlambda} If
  $\bm 1\transpose \caprvec_\userind< \ratemin$, then equation
  \eqref{eq:ubzrveclambda} has no roots. If
  $\bm 1\transpose \caprvec_\userind\geq \ratemin$, then
  \eqref{eq:ubzrveclambda} has a unique root. This root lies in the
  interval
  $[\lambdaslow_\userind\itnot{\itind},
  \lambdashigh_\userind\itnot{\itind}]$, where
  \begin{subequations}
    \begin{align}
      \label{eq:ublambdaslow}
    \lambdaslow_\userind\itnot{\itind}=&
    \min_\gridptind[\raters_\userind\itnot{\itind+1}\entnot{\gridptind} +
    \urs_\userind\itnot{\itind}\entnot{\gridptind} -\caprs_\userind\entnot{\gridptind}]
      \\
      \label{eq:ublambdashigh}
      \lambdashigh_\userind\itnot{\itind}=&
\max\{ \raters_\userind\itnot{\itind+1}\entnot{\gridptind} +
                                           \urs_\userind\itnot{\itind}\entnot{\gridptind}:\gridptind\in                 \{
\gridptind: \caprs_\userind\entnot{\gridptind}>\frac{\ratemin}{\gridptnum}
                                           \}
                                           \}\nonumber\\&-\frac{\ratemin}{\gridptnum}
  \end{align}
\end{subequations}
\end{myproposition}

\balance

\begin{IEEEproof}
  Denote by $\Gfun(\lambdas)$ the left-hand side of
  \eqref{eq:ubzrveclambda}, i.e.,
\begin{align}
  %  \label{eq:ubcombzrs}
   \Gfun(\lambdas) \define \sum_\gridptind
\max(0,    \min(\caprs_\userind\entnot{\gridptind}, \raters_\userind\itnot{\itind+1}\entnot{\gridptind} +
    \urs_\userind\itnot{\itind}\entnot{\gridptind} - \lambdas)).
    \end{align}
    This is a sum of non-increasing piecewise continuous functions and
    therefore $\Gfun$ is also non-increasing piecewise continuous. The
    maximum value is attained for sufficiently small $\lambdas$ and
    equals
    $ \sum_\gridptind\caprs_\userind\entnot{\gridptind} = \bm
    1\transpose \caprvec_\userind$. If
    $ \bm 1\transpose \caprvec_\userind < \ratemin$, then
    $\Gfun(\lambdas)<\ratemin~\forall \lambda$ and
    \eqref{eq:ubzrveclambda} admits no solution. Conversely, if
    $ \bm 1\transpose \caprvec_\userind > \ratemin$, then a solution
    can be found since $\Gfun(\lambdas)>\ratemin$ for sufficiently
    small $\lambdas$ and $\Gfun(\lambdas)=0$ for sufficiently large
    $\lambdas$. Uniqueness follows from the fact that $\Gfun$ is
    strictly decreasing except when $\Gfun(\lambdas)=0$ or
    $\Gfun(\lambdas)=\bm 1\transpose \caprvec_\userind$.
    
      \begin{bullets}%
    \blt[lower]    To show that
    $\Gfun(\lambdaslow_\userind\itnot{\itind})\geq \ratemin$ just note
    from \eqref{eq:ublambdaslow} that
    $\lambdaslow_\userind\itnot{\itind} \leq
    \raters_\userind\itnot{\itind+1}\entnot{\gridptind} +
    \urs_\userind\itnot{\itind}\entnot{\gridptind}
    -\caprs_\userind\entnot{\gridptind}$ or, equivalently,
    $ \caprs_\userind\entnot{\gridptind} \leq
    \raters_\userind\itnot{\itind+1}\entnot{\gridptind} +
    \urs_\userind\itnot{\itind}\entnot{\gridptind}
    -\lambdaslow_\userind\itnot{\itind} $. This clearly yields
    $ \Gfun(\lambdaslow_\userind\itnot{\itind}) = \sum_\gridptind
    \max(0, \caprs_\userind\entnot{\gridptind}) = \sum_\gridptind
    \caprs_\userind\entnot{\gridptind} $, which is greater than or
    equal to $\ratemin$ by assumption.

    \blt[upper]To show that
    $\Gfun(\lambdashigh_\userind\itnot{\itind})\leq \ratemin$, note
    from \eqref{eq:ublambdashigh} that
    $ \lambdashigh_\userind\itnot{\itind} \geq
    \raters_\userind\itnot{\itind+1}\entnot{\gridptind} +
    \urs_\userind\itnot{\itind}\entnot{\gridptind}-{\ratemin}/{\gridptnum}$ for all
    $\gridptind$ such that
    $
    \caprs_\userind\entnot{\gridptind}>{\ratemin}/{\gridptnum}$. This
    clearly implies that
$     \raters_\userind\itnot{\itind+1}\entnot{\gridptind} +
    \urs_\userind\itnot{\itind}\entnot{\gridptind} - \lambdashigh_\userind\itnot{\itind}\leq {\ratemin}/{\gridptnum}$ for all $\gridptind$ such that
    $
    \caprs_\userind\entnot{\gridptind}>{\ratemin}/{\gridptnum}$ and, as a consequence, 
    $ \min(\caprs_\userind\entnot{\gridptind},
    \raters_\userind\itnot{\itind+1}\entnot{\gridptind} +
    \urs_\userind\itnot{\itind}\entnot{\gridptind} - \lambdas))\leq
    {\ratemin}/{\gridptnum}$ and the inequality $\Gfun(\lambdashigh_\userind\itnot{\itind})\leq \ratemin$ follows.

  \end{bullets}%
    
\end{IEEEproof}
  \end{bullets}%

  \blt[U-step]\textbf{$\Umat$-step.} Finally, the $\Umat$-update in \eqref{eq:admmstditu}
  for the assignments in \eqref{eq:ubassignments} becomes
  \begin{align}
    \Umat\itnot{\itind+1} = \Umat\itnot{\itind} +
    \ratemat\itnot{\itind+1}  -
    \Zmat\itnot{\itind+1}.
  \end{align}

\end{bullets}%

%\subsection{Additional Numerical Experiment}


\begin{thebibliography}{10}

\bibitem{zeng2019accessing}
Y.~Zeng, Q.~Wu, and R.~Zhang,
\newblock ``Accessing from the sky: A tutorial on {UAV} communications for {5G}
  and beyond,''
\newblock {\em arXiv preprint arXiv:1903.05289}, 2019.

\bibitem{han2009manet}
Z.~Han, A.~L. Swindlehurst, and K.~J.~R. Liu,
\newblock ``Optimization of manet connectivity via smart deployment/movement of
  unmanned air vehicles,''
\newblock {\em IEEE Trans. Veh. Technol.}, vol. 58, no. 7, pp. 3533--3546,
  2009.

\bibitem{boryaliniz2016placement}
I.~Bor-Yaliniz, A.~El-Keyi, and H.~Yanikomeroglu,
\newblock ``Efficient 3-d placement of an aerial base station in next
  generation cellular networks,''
\newblock in {\em Proc. IEEE Int. Conf. Commun.} IEEE, 2016, pp. 1--5.

\bibitem{chen2017map}
J.~Chen and D.~Gesbert,
\newblock ``Optimal positioning of flying relays for wireless networks: A {LOS}
  map approach,''
\newblock in {\em Proc. IEEE Int. Conf. Commun.}, Paris, France, May 2017, pp.
  1--6.

\bibitem{wang2018adaptive}
Z.~Wang, L.~Duan, and R.~Zhang,
\newblock ``Adaptive deployment for {UAV}-aided communication networks,''
\newblock {\em IEEE Trans. Wireless Commun.}, vol. 18, no. 9, pp. 4531--4543,
  2019.

\bibitem{romero2019noncooperative}
D.~Romero and G.~Leus,
\newblock ``Non-cooperative aerial base station placement via stochastic
  optimization,''
\newblock in {\em Proc. IEEE Mobile Ad-hoc Sensor Netw.}, Shenzhen, China, Dec.
  2019, pp. 131--136.

\bibitem{park2018formation}
S.~Park, K.~Kim, H.~Kim, and H.~Kim,
\newblock ``Formation control algorithm of multi-uav-based network
  infrastructure,''
\newblock {\em Applied Sciences}, vol. 8, no. 10, pp. 1740, 2018.

\bibitem{kim2018topology}
D.-Y. Kim and J.-W. Lee,
\newblock ``Integrated topology management in flying ad hoc networks: Topology
  construction and adjustment,''
\newblock {\em IEEE Access}, vol. 6, pp. 61196--61211, 2018.

\bibitem{alhourani2014urban}
A.~Al-Hourani, S.~Kandeepan, and A.~Jamalipour,
\newblock ``Modeling air-to-ground path loss for low altitude platforms in
  urban environments,''
\newblock in {\em IEEE Global Commun. Conf.}, 2014, pp. 2898--2904.

\bibitem{kalantari2016number}
E.~Kalantari, H.~Yanikomeroglu, and A.~Yongacoglu,
\newblock ``On the number and {3D} placement of drone base stations in wireless
  cellular networks,''
\newblock in {\em IEEE Vehicular Tech. Conf.}, 2016, pp. 1--6.

\bibitem{hammouti2019mechanism}
H.~El~Hammouti, M.~Benjillali, B.~Shihada, and M.-S. Alouini,
\newblock ``A distributed mechanism for joint {3D} placement and user
  association in {UAV}-assisted networks,''
\newblock in {\em IEEE Wireless Commun. Netw. Conf.}, Marrakech, Morocco, Apr.
  2019.

\bibitem{perabathini2019qos}
B.~Perabathini, K.~Tummuri, A.~Agrawal, and V.S. Varma,
\newblock ``Efficient {3D} placement of {UAVs with QoS Assurance in Ad Hoc
  Wireless Networks},''
\newblock in {\em Int. Conf. Comput. Commun. Netw.}, 2019, pp. 1--6.

\bibitem{liu2019deployment}
X.~Liu, Y.~Liu, and Y.~Chen,
\newblock ``Reinforcement learning in multiple-{UAV} networks: Deployment and
  movement design,''
\newblock {\em IEEE Trans. Veh. Tech.}, vol. 68, no. 8, pp. 8036--8049, 2019.

\bibitem{shehzad2021backhaul}
M.K. Shehzad, A.~Ahmad, S.A. Hassan, and H.~Jung,
\newblock ``Backhaul-aware intelligent positioning of {UAVs} and association of
  terrestrial base stations for fronthaul connectivity,''
\newblock {\em IEEE Trans. Netw. Sci. Eng.}, pp. 1--1, 2021.

\bibitem{qiu2020reinforcement}
J.~Qiu, J.~Lyu, and L.~Fu,
\newblock ``Placement optimization of aerial base stations with deep
  reinforcement learning,''
\newblock in {\em IEEE Int. Conf. Commun.}, 2020, pp. 1--6.

\bibitem{sabzehali2021orientation}
J.~Sabzehali, V.K. Shah, H.S. Dhillon, and J.H. Reed,
\newblock ``{3D} placement and orientation of {mmWave-based UAVs for Guaranteed
  LoS Coverage},''
\newblock {\em IEEE Wireless Commun. Letters}, pp. 1--1, 2021.

\bibitem{patwari2008nesh}
N.~Patwari and P.~Agrawal,
\newblock ``Nesh: A joint shadowing model for links in a multi-hop network,''
\newblock in {\em Proc. IEEE Int. Conf. Acoust., Speech, Signal Process.}, Las
  Vegas, NV, Mar. 2008, pp. 2873--2876.

\bibitem{patwari2008correlated}
N.~Patwari and P.~Agrawal,
\newblock ``Effects of correlated shadowing: Connectivity, localization, and
  {RF} tomography,''
\newblock in {\em Proc. Int. Conf. Info. Process. Sensor Networks}, St. Louis,
  MO, Apr. 2008, pp. 82--93.

\bibitem{romero2018blind}
D.~Romero, D.~Lee, and G.~B. Giannakis,
\newblock ``Blind radio tomography,''
\newblock {\em IEEE Trans. Signal Process.}, vol. 66, no. 8, pp. 2055--2069,
  Jan. 2018.

\bibitem{wilson2009regularization}
J.~Wilson, N.~Patwari, and O.~G. Vasquez,
\newblock ``Regularization methods for radio tomographic imaging,''
\newblock in {\em Virginia Tech Symp. Wireless Personal Commun.}, Blacksburg,
  VA, Jun. 2009.

\bibitem{kanso2009compressed}
M.~A. Kanso and M.~G. Rabbat,
\newblock ``Compressed rf tomography for wireless sensor networks: Centralized
  and decentralized approaches,''
\newblock in {\em Int. Conf. Distributed Comput. Sensor Syst.}, Marina del Rey,
  CA, 2009, Springer, pp. 173--186.

\bibitem{hamilton2014modeling}
B.~R. Hamilton, X.~Ma, R.~J. Baxley, and S.~M. Matechik,
\newblock ``Propagation modeling for radio frequency tomography in wireless
  networks,''
\newblock {\em IEEE J. Sel. Topics Signal Process.}, vol. 8, no. 1, pp. 55--65,
  Feb. 2014.

\bibitem{mitchell1990comparison}
J.R. Mitchell, P.~Dickof, and A.G. Law,
\newblock ``A comparison of line integral algorithms,''
\newblock {\em Comput. Physics}, vol. 4, no. 2, pp. 166--172, 1990.

\bibitem{lin2021admm}
T.~Lin, S.~Ma, Y.~Ye, and S.~Zhang,
\newblock ``An {ADMM}-based interior-point method for large-scale linear
  programming,''
\newblock {\em Optim. Methods Software}, vol. 36, no. 2-3, pp. 389--424, 2021.

\bibitem{boyd2011distributed}
S.~Boyd, N.~Parikh, E.~Chu, B.~Peleato, and J.~Eckstein,
\newblock ``Distributed optimization and statistical learning via the
  alternating direction method of multipliers,''
\newblock {\em Found. Trends Mach. Learn.}, vol. 3, no. 1, pp. 1--122, Jan.
  2011.

\bibitem{candes2008reweighted}
E.J. Candes, M.B. Wakin, and S.P. Boyd,
\newblock ``Enhancing sparsity by reweighted $\ell_1$ minimization,''
\newblock {\em J. Fourier Analysis App.}, vol. 14, no. 5, pp. 877--905, 2008.

\bibitem{huang2020sparse}
M.~Huang, L.~Huang, S.~Zhong, and P.~Zhang,
\newblock ``{UAV}-mounted mobile base station placement via sparse recovery,''
\newblock {\em IEEE Access}, vol. 8, pp. 71775--71781, 2020.

\bibitem{galkin2016deployment}
B.~Galkin, J.~Kibilda, and L.A. DaSilva,
\newblock ``Deployment of {UAV}-mounted access points according to spatial user
  locations in two-tier cellular networks,''
\newblock in {\em Wireless Days}. IEEE, 2016, pp. 1--6.

\bibitem{lyu2017mounted}
J.~Lyu, Y.~Zeng, R.~Zhang, and T.J. Lim,
\newblock ``Placement optimization of {UAV}-mounted mobile base stations,''
\newblock {\em IEEE Commun. Letters}, vol. 21, no. 3, pp. 604--607, 2017.

\bibitem{chen2021relay}
J.~Chen, U.~Mitra, and D.~Gesbert,
\newblock ``{3D} urban {UAV} relay placement: Linear complexity algorithm and
  analysis,''
\newblock {\em IEEE Trans. Wireless Commun.}, pp. 1--1, 2021.

\end{thebibliography}
\end{document}